\documentclass[10pt]{article}
\usepackage{setspace}
 
\usepackage{amsmath} 
\usepackage{amssymb}
\usepackage{latexsym}
\usepackage{xcolor}
\usepackage{fancyhdr}
\allowdisplaybreaks 

 \usepackage{comment}

\def\XXint#1#2#3{{\setbox0=\hbox{$#1{#2#3}{\int}$} 
\vcenter{\hbox{$#2#3$}}\kern-.5\wd0}}   

 \numberwithin{equation}{section}
\newtheorem{theorem}[equation]{Theorem}
\newtheorem{proposition}[equation]{Proposition}
\newtheorem{definition}[equation]{Definition}
\newtheorem{remark}[equation]{Remark}
\newtheorem{lemma}[equation]{Lemma}

\title{
The volume potential for elliptic differential operators in 
Schauder spaces } 
 
\author{  
Massimo Lanza de Cristoforis
\\
Dipartimento di Matematica `Tullio Levi-Civita', 
\\
Universit\`a degli Studi di Padova, 
\\
Via Trieste 63, Padova 35121, 
Italy. 
\\
E-mail: mldc@math.unipd.it   }

\date{\ }

\begin{document}
 
 \maketitle

\noindent
{\bf Abstract:}  The aim of this paper is to prove continuity results for 
  the volume potential corresponding to the fundamental solution of a  second order differential operator with constant coefficients in Schauder spaces   of   negative exponent and to generalize some classical results
  in  Schauder spaces of positive exponents.

 \vspace{\baselineskip}



\section{Introduction} The aim of this paper is to prove continuity results for 
  the volume potential corresponding to the fundamental solution of a  second order differential operator with constant coefficients in Schauder spaces   of   negative exponent and to generalize some classical results
  in  Schauder spaces of positive exponents. Unless otherwise specified,  we assume  throughout the paper that
\[
n\in {\mathbb{N}}\setminus\{0,1\}\,,
\]
where ${\mathbb{N}}$ denotes the set of natural numbers including $0$. Let $\alpha\in[0,1]$, $m\in {\mathbb{N}}\setminus\{0\}$. Let $\Omega$ be a bounded open subset of ${\mathbb{R}}^{n}$ of class $C^{m,\alpha}$. For the definition and properties of the classical  Schauder spaces both of negative and positive exponent, we refer for example to \cite[Chap.~2]{DaLaMu21}. We also find convenient to set 
\[
\Omega^-\equiv {\mathbb{R}}^n\setminus\overline{\Omega}\,,
\] 
where $\overline{\Omega}$ denotes the closure of $\Omega$. We employ the same notation of reference \cite{DoLa17} with Dondi that we now introduce. 
 Let $N_{2}$ denote the number of multi-indexes $\gamma\in {\mathbb{N}}^{n}$ with $|\gamma|\leq 2$. For each 
\begin{equation}
\label{introd0}
{\mathbf{a}}\equiv (a_{\gamma})_{|\gamma|\leq 2}\in {\mathbb{C}}^{N_{2}}\,, 
\end{equation}
we set 
\[
a^{(2)}\equiv (a_{lj} )_{l,j=1,\dots,n}\qquad
a^{(1)}\equiv (a_{j})_{j=1,\dots,n}\qquad
a\equiv a_{0} 
\]
with $a_{lj} \equiv 2^{-1}a_{e_{l}+e_{j}}$ for $j\neq l$, $a_{jj} \equiv
 a_{e_{j}+e_{j}}$,
and $a_{j}\equiv a_{e_{j}}$, where $\{e_{j}:\,j=1,\dots,n\}$  is the canonical basis of ${\mathbb{R}}^{n}$. We note that the matrix $a^{(2)}$ is symmetric. 
Then we assume that 
  ${\mathbf{a}}\in  {\mathbb{C}}^{N_{2}}$ satisfies the following ellipticity assumption
\begin{equation}
\label{ellip}
\inf_{
\xi\in {\mathbb{R}}^{n}, |\xi|=1
}{\mathrm{Re}}\,\left\{
 \sum_{|\gamma|=2}a_{\gamma}\xi^{\gamma}\right\} >0\,,
\end{equation}
and we consider  the case in which
\begin{equation}
\label{symr}
a_{lj} \in {\mathbb{R}}\qquad\forall  l,j=1,\dots,n\,.
\end{equation}
Then we introduce the differential operator 
\begin{eqnarray*}
P[{\mathbf{a}},D]u&\equiv&\sum_{l,j=1}^{n}\partial_{x_{l}}(a_{lj}\partial_{x_{j}}u)
+
\sum_{l=1}^{n}a_{l}\partial_{x_{l}}u+au\,,
\end{eqnarray*}
for all $u,v\in C^{2}(\overline{\Omega})$, and a fundamental solution $S_{{\mathbf{a}} }$ of $P[{\mathbf{a}},D]$. Then we can consider the restriction map $r_{|\overline{\Omega}}$ from the space of test functions ${\mathcal{D}}({\mathbb{R}}^n)$ to the Schauder space  $C^{1,\alpha}(\overline{\Omega})$. Then the transpose map $r_{|\overline{\Omega}}^t$ is linear and continuous from $(C^{1,\alpha}(\overline{\Omega}))'$ to ${\mathcal{D}}'({\mathbb{R}}^n)$. Moreover, if $\mu\in (C^{1,\alpha}(\overline{\Omega}))'$, then $r_{|\overline{\Omega}}^t\mu$ has compact support. Hence, it makes sense to consider the convolution of 
  $r_{|\overline{\Omega}}^t\mu$ with  the fundamental solution $S_{{\mathbf{a}} }$. Namely, the   distribution
\[
{\mathcal{P}}_\Omega[S_{{\mathbf{a}} },\mu]=(r_{|\overline{\Omega}}^t\mu)\ast S_{{\mathbf{a}} }\in {\mathcal{D}}'({\mathbb{R}}^n)\,.
\]
Then we    set
\begin{eqnarray}\label{prop:dvpsa3}
{\mathcal{P}}_\Omega^+[S_{{\mathbf{a}} },\mu]&\equiv&\left((r_{|\overline{\Omega}}^t\mu)\ast S_{{\mathbf{a}} }\right)_{|\Omega}
\qquad\text{in}\ \Omega\,,
\\ \nonumber
{\mathcal{P}}_\Omega^-[S_{{\mathbf{a}} },\mu] &\equiv&
\left((r_{|\overline{\Omega}}^t\mu)\ast S_{{\mathbf{a}} }\right)_{|\Omega^-}
\qquad\text{in}\ \Omega^-\,.
\end{eqnarray}
Next we  consider the classical volume potential and we generalize to nonhomogeneous differential operators  as $P[{\mathbf{a}},D]$ and to case $\alpha=1$, a known result of Miranda \cite[Thm.~3.I, p.~320]{Mi65} for homogeneous differential operators in case $\alpha\in]0,1[$. See also Kirsch and Hettlich \cite[\S 3.1.2]{KiHe15} for the Helmholtz operator.\par 
  
   Namely we take $\mu
 \in C^{m,\alpha}(\overline{\Omega})$ with $m\in {\mathbb{N}}$, $\alpha\in]0,1]$, we associate to $\mu$ an  element ${\mathcal{J}}[\mu]$ of $(C^{1,\alpha}(\overline{\Omega}))'$ (cf.~Lemma \ref{lem:cainclo}) and we prove that if $\Omega$ is bounded and of class $C^{m+1,\alpha}$, then 
  the map ${\mathcal{P}}_\Omega^+[S_{{\mathbf{a}} },{\mathcal{J}}[\cdot]]$ is linear and continuous from $C^{m,\alpha}(\overline{\Omega}) $ to $C^{m+2,\alpha}(\overline{\Omega})$ in case $\alpha\in]0,1[$  and to the generalized Schauder space
$C^{m+2,\omega_{1}(\cdot)}(\overline{\Omega})$ of functions with $(m+2)$-th order derivatives which satisfy a generalized $\omega_{1}(\cdot)$-H\"{o}lder condition with $\omega_{1}(\cdot)$ as in (\ref{omth}) below and thus with 
\[
\omega_{1}(r)\sim r|\ln r| 
\qquad{\mathrm{as}}\ r\to 0,
\]
 in case $\alpha=1$. Similarly, we prove that 
    if $r\in ]0,+\infty[$ is such that $\overline{\Omega}\subseteq {\mathbb{B}}_n(0,r)$, then the  map ${\mathcal{P}}_\Omega^-[S_{{\mathbf{a}} },{\mathcal{J}}[\cdot]]_{|\overline{{\mathbb{B}}_n(0,r)}\setminus\Omega}$ is linear and continuous from the space $C^{m,\alpha}(\overline{\Omega}) $ to $C^{m+2,\alpha}(\overline{{\mathbb{B}}_n(0,r)}\setminus\Omega)$ in case $\alpha\in]0,1[$ and to $C^{m+2,\omega_1(\cdot)}(\overline{{\mathbb{B}}_n(0,r)}\setminus\Omega)$ in case $\alpha=1$ (see Theorem \ref{thm:vopoadma}).\par
   
Next we turn to consider the Schauder space $C^{-1,\alpha}(\overline{\Omega})$ of sums of partial distributional  derivatives of order up to one of $\alpha$-H\"{o}lder continuous functions in $\Omega$ (cf.~\textit{e.g.}, Dalla Riva, the author and Musolino  \cite[\S 2.22]{DaLaMu21}). 

Here we mention that the space $C^{-1,\alpha}(\overline{\Omega})$  has been known for a long time and has been used in the analysis of elliptic and parabolic partial differential equations (cf. Triebel \cite{Tr78}, Gilbarg and Trudinger~\cite{GiTr83}, 
Vespri \cite{Ve88}, Lunardi and Vespri \cite{LuVe91}, Dalla Riva, the author and Musolino \cite{DaLaMu21}, \cite{La08a}).

One can prove that there exists an extension operator $E^\sharp$ from $C^{-1,\alpha}(\overline{\Omega})$ to 
  $(C^{1,\alpha}(\overline{\Omega}))'$ (see \cite[\S 3]{La24c}).\par
  
Then we prove that if $\Omega$ is bounded and of class $C^{1,\alpha}$ with $\alpha\in]0,1[$, then the map ${\mathcal{P}}_\Omega^+[S_{{\mathbf{a}} },E^\sharp[\cdot]]$ is linear and continuous from $C^{-1,\alpha}(\overline{\Omega}) $ to $C^{1,\alpha}(\overline{\Omega})$ and that if $\Omega$ is a bounded Lipschitz subset of ${\mathbb{R}}^n$, then  the map ${\mathcal{P}}_\Omega^+[S_{{\mathbf{a}} },E^\sharp[\cdot]]$ is linear and continuous from $C^{-1,1}(\overline{\Omega}) $ to
  to $C^{1,\omega_1(\cdot)}(\overline{\Omega})$. Similarly, we prove a corresponding statement for  ${\mathcal{P}}_\Omega^-[E^\sharp[\cdot]]_{|\overline{{\mathbb{B}}_n(0,r)}\setminus\Omega}$,  where $r\in ]0,+\infty[$ is such that $\overline{\Omega}\subseteq {\mathbb{B}}_n(0,r)$ (see Propositions \ref{prop:dvpsnecr-1a}, \ref{prop:dvpsnecr-11}). Such results   extend    a corresponding result of \cite[Thm.~3.6 (ii)]{La08a}, Dalla Riva, the author and Musolino 
  \cite[Thm.~7.19]{DaLaMu21} for the Laplace operator in case $\alpha\in]0,1[$ and find  application in the analysis of a nonvariational form of the Neumann problem for the Poisson equation (see \cite[\S 6]{La24c}). We also mention the extension to the case of the heat volume potential of Luzzini \cite{Lu24}.

The paper is organized as follows. Section \ref{sec:tecprel}   is a section of preliminaries and notation. In Section \ref{sec:prelfuso}, we introduce some properties on the fundamental solution $S_{{\mathbf{a}} }$ that we need. In Section \ref{sec:dvopo} we collect some preliminary properties of the distributional volume potential. 
In Section \ref{sec:teleinop} we prove a variant of a technical statement on an integral operator that has been proved in   Dalla Riva, the author  and Musolino \cite[Prop.~7.15]{DaLaMu21}.
In Section \ref{sec:voscpoe} we prove our generalization of the result of Miranda \cite[Thm.~3.I, p.~320]{Mi65} for Schauder spaces with positive exponents. In Section \ref{sec:voscnege} we prove our results in case the density belongs to a Schauder space with  negative  exponent. In the appendix at the end of the paper, we prove a formula of integration by parts for kernels with a weak singularity. Related formulas
 are known even in case of sets with a rough boundary. See for example Mitrea, Mitrea and Mitrea \cite[Thm.~1.11.8]{MitMitMit22}. Then we also include an   extension of Dalla Riva, the author and Musolino \cite{DaLaMu24a} of a  result  of  Miranda \cite{Mi65}    for singular integrals and   an extension of \cite{La24e} of a   result for single layer potentials of  Miranda \cite{Mi65}.
 
 \section{Preliminaries and notation}\label{sec:tecprel}
If $X$ and $Y$, $Z$ are normed spaces, then ${\mathcal{L}}(X,Y)$ denotes the space of linear and continuous maps from $X$ to $Y$ and ${\mathcal{L}}^{(2)}(X\times Y, Z)$ 
denotes the space of bilinear and continuous maps from $X\times Y$ to $Z$ with their usual operator norm (cf.~\textit{e.g.}, \cite[pp.~16, 621]{DaLaMu21}). 
$|A|$ denotes the operator norm of a matrix $A$ with real (or complex) entries, 
       $A^{t}$ denotes the transpose matrix of $A$. $\delta_{l,j}$ denotes the Kronecker\index{Kronecker symbol}  symbol. Namely,  $\delta_{l,j}=1$ if $l=j$, $\delta_{l,j}=0$ if $l\neq j$, with $l,j\in {\mathbb{N}}$.  The symbol
$| \cdot|$ denotes the Euclidean modulus   in
${\mathbb{R}}^{n}$ or in ${\mathbb{C}}$. For all $r\in]0,+\infty[$, $ x\in{\mathbb{R}}^{n}$, 
$x_{j}$ denotes the $j$-th coordinate of $x$, and  
 ${\mathbb{B}}_{n}( x,r)$ denotes the ball $\{y\in{\mathbb{R}}^{n}:\, |x- y|<r\}$.  If ${\mathbb{D}}$ is a subset of $ {\mathbb{R}}^n$, 
then we set
\[
B({\mathbb{D}})\equiv\left\{
f\in {\mathbb{C}}^{\mathbb{D}}:\,f\ \text{is\ bounded}
\right\}
\,,\quad
\|f\|_{B({\mathbb{D}})}\equiv\sup_{\mathbb{D}}|f|\qquad\forall f\in B({\mathbb{D}})\,.
\]
Then $C^0({\mathbb{D}})$ denotes the set of continuous functions from ${\mathbb{D}}$ to ${\mathbb{C}}$ and we introduce the subspace
$
C^0_b({\mathbb{D}})\equiv C^0({\mathbb{D}})\cap B({\mathbb{D}})
$
of $B({\mathbb{D}})$.  Let $\omega$ be a function from $[0,+\infty[$ to itself such that
\begin{eqnarray}
\nonumber
&&\qquad\qquad\omega(0)=0,\qquad \omega(r)>0\qquad\forall r\in]0,+\infty[\,,
\\
\label{om}
&&\qquad\qquad\omega\ {\text{is\   increasing,}}\ \lim_{r\to 0^{+}}\omega(r)=0\,,
\\
\nonumber
&&\qquad\qquad{\text{and}}\ \sup_{(a,t)\in[1,+\infty[\times]0,+\infty[}
\frac{\omega(at)}{a\omega(t)}<+\infty\,.
\end{eqnarray}
Here `$\omega$ is increasing' means that 
$\omega(r_1)\leq \omega(r_2)$ whenever $r_1$, $r_2\in [0,+\infty[$ and $r_1<r_2$.
If $f$ is a function from a subset ${\mathbb{D}}$ of ${\mathbb{R}}^n$   to ${\mathbb{C}}$,  then we denote by   $|f:{\mathbb{D}}|_{\omega (\cdot)}$  the $\omega(\cdot)$-H\"older constant  of $f$, which is delivered by the formula   
\[
|f:{\mathbb{D}}|_{\omega (\cdot)
}
\equiv
\sup\left\{
\frac{|f( x )-f( y)|}{\omega(| x- y|)
}: x, y\in {\mathbb{D}} ,  x\neq
 y\right\}\,.
\]        
If $|f:{\mathbb{D}}|_{\omega(\cdot)}<+\infty$, we say that $f$ is $\omega(\cdot)$-H\"{o}lder continuous. Sometimes, we simply write $|f|_{\omega(\cdot)}$  
instead of $|f:{\mathbb{D}}|_{\omega(\cdot)}$. The
subset of $C^{0}({\mathbb{D}} ) $  whose
functions  are
$\omega(\cdot)$-H\"{o}lder continuous    is denoted  by  $C^{0,\omega(\cdot)} ({\mathbb{D}})$
and $|f:{\mathbb{D}}|_{\omega(\cdot)}$ is a semi-norm on $C^{0,\omega(\cdot)} ({\mathbb{D}})$.  
Then we consider the space  $C^{0,\omega(\cdot)}_{b}({\mathbb{D}} ) \equiv C^{0,\omega(\cdot)} ({\mathbb{D}} )\cap B({\mathbb{D}} ) $ with the norm \[
\|f\|_{ C^{0,\omega(\cdot)}_{b}({\mathbb{D}} ) }\equiv \sup_{x\in {\mathbb{D}} }|f(x)|+|f|_{\omega(\cdot)}\qquad\forall f\in C^{0,\omega(\cdot)}_{b}({\mathbb{D}} )\,.
\] 
\begin{remark}
\label{rem:om4}
Let $\omega$ be as in (\ref{om}). 
Let ${\mathbb{D}}$ be a   subset of ${\mathbb{R}}^{n}$. Let $f$ be a bounded function from $ {\mathbb{D}}$ to ${\mathbb{C}}$, $a\in]0,+\infty[$.  Then,
\[
\label{rem:om5}
\sup_{x,y\in {\mathbb{D}},\ |x-y|\geq a}\frac{|f(x)-f(y)|}{\omega(|x-y|)}
\leq \frac{2}{\omega(a)} \sup_{{\mathbb{D}}}|f|\,.
\]
\end{remark}
In the case in which $\omega(\cdot)$ is the function 
$r^{\alpha}$ for some fixed $\alpha\in]0,1]$, a so-called H\"{o}lder exponent, we simply write $|\cdot:{\mathbb{D}}|_{\alpha}$ instead of
$|\cdot:{\mathbb{D}}|_{r^{\alpha}}$, $C^{0,\alpha} ({\mathbb{D}})$ instead of $C^{0,r^{\alpha}} ({\mathbb{D}})$, $C^{0,\alpha}_{b}({\mathbb{D}})$ instead of $C^{0,r^{\alpha}}_{b} ({\mathbb{D}})$, and we say that $f$ is $\alpha$-H\"{o}lder continuous provided that 
$|f:{\mathbb{D}}|_{\alpha}<+\infty$. For each $\theta\in]0,1]$, we define the function $\omega_{\theta}(\cdot)$ from $[0,+\infty[$ to itself by setting
\begin{equation}
\label{omth}
\omega_{\theta}(r)\equiv
\left\{
\begin{array}{ll}
0 &r=0\,,
\\
r^{\theta}|\ln r | &r\in]0,r_{\theta}]\,,
\\
r_{\theta}^{\theta}|\ln r_{\theta} | & r\in ]r_{\theta},+\infty[\,,
\end{array}
\right.
\end{equation}
where
$
r_{\theta}\equiv e^{-1/\theta}
$ for all $\theta\in ]0,1]$. Obviously, $\omega_{\theta}(\cdot) $ is concave and satisfies   condition (\ref{om}).
 We also note that if ${\mathbb{D}}\subseteq {\mathbb{R}}^n$, then the continuous embeddings
\[
C^{0, \theta }_b({\mathbb{D}})\subseteq 
C^{0,\omega_\theta(\cdot)}_b({\mathbb{D}})\subseteq 
C^{0,\theta'}_b({\mathbb{D}})
\]
hold  for all $\theta'\in ]0,\theta[$. For the standard properties of the spaces of H\"{o}lder or Lipschitz continuous functions, we refer    to \cite[\S 2]{DoLa17}, \cite[\S 2.6]{DaLaMu21}.\par

Let $\Omega$ be an open subset of ${\mathbb{R}}^n$. The space of $m$ times continuously 
differentiable complex-valued functions on $\Omega$ is denoted by 
$C^{m}(\Omega,{\mathbb{C}})$, or more simply by $C^{m}(\Omega)$. 
Let $f\in  C^{m}(\Omega)  $. Then   $Df$ denotes the Jacobian matrix of $f$. 
Let  $\eta\equiv
(\eta_{1},\dots ,\eta_{n})\in{\mathbb{N}}^{n}$, $|\eta |\equiv
\eta_{1}+\dots +\eta_{n}  $. Then $D^{\eta} f$ denotes
$\frac{\partial^{|\eta|}f}{\partial
x_{1}^{\eta_{1}}\dots\partial x_{n}^{\eta_{n}}}$.    The
subspace of $C^{m}(\Omega )$ of those functions $f$ whose derivatives $D^{\eta }f$ of
order $|\eta |\leq m$ can be extended with continuity to 
$\overline{\Omega}$  is  denoted $C^{m}(
\overline{\Omega})$. 
The
subspace of $C^{m}(\overline{\Omega} ) $  whose
functions have $m$-th order derivatives that are
H\"{o}lder continuous  with exponent  $\alpha\in
]0,1]$ is denoted $C^{m,\alpha} (\overline{\Omega})$ and the
subspace of $C^{m}(\overline{\Omega} ) $  whose
functions have $m$-th order derivatives that are
$\omega(\cdot)$-H\"{o}lder continuous    is denoted $C^{m,\omega(\cdot)} (\overline{\Omega})$. 

The subspace of $C^{m}(\overline{\Omega} ) $ of those functions $f$ such that the restriction $f_{|\overline{(\Omega\cap{\mathbb{B}}_{n}(0,r))}}$ belongs to $
C^{m,\omega(\cdot)}(\overline{(\Omega\cap{\mathbb{B}}_{n}(0,r))})$ with $\omega$ as in (\ref{om}) 
for all $r\in]0,+\infty[$ is denoted $C^{m,\omega(\cdot)}_{{\mathrm{loc}}}(\overline{\Omega} ) $. \par
Now let $\Omega $ be a bounded
open subset of  ${\mathbb{R}}^{n}$. Then $C^{m}(\overline{\Omega} )$,  $C^{m,\omega(\cdot) }(\overline{\Omega})$ with $\omega$ as in (\ref{om}),
  $C^{m,\alpha }(\overline{\Omega})$,  are endowed with their usual norm and are well known to be 
Banach spaces  (cf.~\textit{e.g.},  \cite[\S 2]{DoLa17}, Dalla Riva, the author and Musolino  \cite[\S 2.11]{DaLaMu21}).  

 For the definition of a bounded open Lipschitz subset of ${\mathbb{R}}^{n}$ and  for  the (classical) definition of open set of class $C^{m}$ or of class $C^{m,\alpha}$, we refer for example to 
 Dalla Riva, the author and Musolino  \cite[\S 2.9, \S 2.13]{DaLaMu21}. 
 
  For the (classical) definition of the (generalized) H\"{o}lder and Schauder spaces  $C^{m,\omega(\cdot) }(\partial\Omega)$ with $\omega$ as in (\ref{om}),  $C^{m,\alpha}(\partial\Omega)$ 
on the boundary $\partial\Omega$ of an open set $\Omega$ for some $m\in{\mathbb{N}}$, $\alpha\in ]0,1]$, we refer for example to Dondi and the  author \cite[\S 2]{DoLa17}, 
    Dalla Riva, the author and Musolino  \cite[\S  2.20]{DaLaMu21}.

 The space of real valued functions of class $C^\infty$ with compact support in an open set $\Omega$ of ${\mathbb{R}}^n$ is denoted ${\mathcal{D}}(\Omega)$. Then its dual ${\mathcal{D}}'(\Omega)$ is known to be the space of distributions in $\Omega$. The support of a function or of a distribution is denoted by the abbreviation `${\mathrm{supp}}$'.  

Morever, we retain the standard notation for the Lebesgue spaces $L^p$ for $p\in [1,+\infty]$ (cf.~\textit{e.g.}, Folland \cite[Chapt.~6]{Fo95}, \cite[\S 2.1]{DaLaMu21}) and
$m_n$ denotes the
$n$ dimensional Lebesgue measure.\par 

If $\Omega$ is a  bounded open subset of ${\mathbb{R}}^{n}$, then we find convenient to consider the dual $(C^{1,\alpha}(\overline{\Omega}))'$  of $C^{1,\alpha}(\overline{\Omega})$ with its usual (normable) topology and the corresponding duality pairing $\langle\cdot,\cdot\rangle$ and we say that the elements of $(C^{1,\alpha}(\overline{\Omega}))'$ are distributions in 
 $\overline{\Omega}$. The following Lemma is well known and is an immediate consequence of the H\"{o}lder inequality.
\begin{lemma}\label{lem:cainclo}
 Let  $\alpha\in ]0,1]$. Let $\Omega$ be a bounded open Lipschitz subset of ${\mathbb{R}}^{n}$.  Then the canonical inclusion  ${\mathcal{J}}$ from the Lebesgue space $L^1(\Omega)$ of integrable functions in $\Omega$ to $(C^{1,\alpha}(\overline{\Omega}))'$ that takes $f$ to the functional ${\mathcal{J}}[f]$ defined by 
 \begin{equation}\label{lem:cainclo1}
\langle{\mathcal{J}}[f],v\rangle \equiv \int_{\Omega}f v\,d\sigma\qquad\forall v\in C^{1,\alpha}(\overline{\Omega})\,,
\end{equation}
is linear continuous and injective.
\end{lemma}
As customary, we say  that ${\mathcal{J}}[f]$ is the `distribution that is canonically associated to $f$' and we  omit the indication of the inclusion map ${\mathcal{J}}$ when no ambiguity can arise. By Lemma \ref{lem:cainclo}, the space $C^{0,\alpha}(\overline{\Omega})$ is continuously embedded into $(C^{1,\alpha}(\overline{\Omega}))'$.\par

We now summarize the definition and some elementary properties of the Schau\-der space $C^{-1,\alpha}(\overline{\Omega})$ 
by following the presentation of Dalla Riva, the author and Musolino \cite[\S 2.22]{DaLaMu21}.
\begin{definition} 
\label{defn:sch-1}\index{Schauder space!with negative exponent}
Let $n\in {\mathbb{N}}\setminus\{0\}$. Let $\alpha\in]0,1]$. Let $\Omega$ be a bounded open subset of ${\mathbb{R}}^{n}$. We denote by $C^{-1,\alpha}(\overline{\Omega})$ the subspace 
 \[
 \left\{
 f_{0}+\sum_{j=1}^{n}\frac{\partial}{\partial x_{j}}f_{j}:\,f_{j}\in 
 C^{0,\alpha}(\overline{\Omega})\ \forall j\in\{0,\dots,n\}
 \right\}\,,
 \]
 of the space of distributions ${\mathcal{D}}'(\Omega)$  in $\Omega$.
\end{definition}
According to the above definition, the space $C^{-1,\alpha}(\overline{\Omega})$ is the image of the linear and continuous map 
\begin{equation}\label{eq:xi0a}
\Xi:\,(C^{0,\alpha}(\overline{\Omega}))^{n+1}\to {\mathcal{D}}'(\Omega)
\end{equation}
that takes an $(n+1)$-tuple $(f_{0},\dots,f_{n})$ to $ f_{0}+\sum_{j=1}^{n}\frac{\partial}{\partial x_{j}}f_{j}$.
Then we set
\begin{eqnarray}
\label{defn:sch-2}
\lefteqn{
\|f\|_{  C^{-1,\alpha}(\overline{\Omega})  }
\equiv\inf\biggl\{\biggr.
\sum_{j=0}^{n}\|f_{j}\|_{ C^{0,\alpha}(\overline{\Omega})  }
:\,
}
\\ \nonumber
&&\qquad\qquad\qquad\qquad
f=f_{0}+\sum_{j=1}^{n}\frac{\partial}{\partial x_{j}}f_{j}\,,\ 
f_{j}\in C^{0,\alpha}(\overline{\Omega})\ \forall j\in \{0,\dots,n\}
\biggl.\biggr\}
\end{eqnarray}
and $(C^{-1,\alpha}(\overline{\Omega}), \|\cdot\|_{  C^{-1,\alpha}(\overline{\Omega})  })$ is known to be a Banach space. Also, the definition of the norm $\|\cdot\|_{  C^{-1,\alpha}(\overline{\Omega})  }$ implies that $C^{0,\alpha}(\overline{\Omega})$ is continuously embedded into $C^{-1,\alpha}(\overline{\Omega})$ and that the partial derivation $\frac{\partial}{\partial x_{j}}$ is continuous from 
$C^{0,\alpha}(\overline{\Omega})$ to $C^{-1,\alpha}(\overline{\Omega})$ for all $j\in\{1,\dots,n\}$. 
If $\alpha\in]0,1[$, then the elements of  $C^{-1,\alpha}(\overline{\Omega})$ are distributions that are not necessarily associated to functions. However, if $\alpha=1$, the elements of $C^{-1,1}(\overline{\Omega})$ are associated to essentially bounded functions in the sense of the following statement. 
\begin{proposition}\label{prop:c-1ainfty}
Let $n\in {\mathbb{N}}\setminus\{0\}$.    Let $\Omega$ be a bounded open  subset of ${\mathbb{R}}^{n}$. Then 
 $C^{-1,1}(\overline{\Omega})$ is continuously embedded into $L^\infty(\Omega)$. 
\end{proposition}
{\bf Proof.} If $f\in C^{-1,1}(\overline{\Omega})$, then there exist $f_{j}\in C^{0,1}(\overline{\Omega})$ for all $j\in \{0,\dots,n\}$ such that
\begin{equation}\label{prop:c-1ainfty1}
f=f_{0}+\sum_{j=1}^{n}\frac{\partial}{\partial x_{j}}f_{j}\,.
\end{equation}
By the Rademacher Theorem, $f_j$ is differentiable almost everywhere in $\Omega$, the distributional derivative $\frac{\partial}{\partial x_{j}}f_{j}$ coincides with the classical almost everywhere defined $j$-th partial derivative of $f_{j}$ and
\[
\left\|\frac{\partial}{\partial x_{j}}f_{j}\right\|_{L^\infty(\Omega)}\leq |f_j:\Omega|_1\leq \|f_j\|_{C^{0,1}(\overline{\Omega})}
\]
for all $j\in \{0,\dots,n\}$. Then the triangular inequality implies that
\[
\|f\|_{L^\infty(\Omega)}\leq \|f_0\|_{L^\infty(\Omega)}+\sum_{j=1}^{n}\left\|\frac{\partial}{\partial x_{j}}f_{j}\right\|_{L^\infty(\Omega)}\leq \|f_0\|_{L^\infty(\Omega)}+\sum_{j=1}^{n}\|f_j\|_{C^{0,1}(\overline{\Omega})}\,.
\]
Then by taking the infimum on all possible $\{f_j\}_{j\in \{0,\dots,n\}}$ as in (\ref{prop:c-1ainfty1}), we deduce that
$\|f\|_{L^\infty(\Omega)}\leq \|f\|_{C^{-1,1}(\overline{\Omega})}$ and thus the proof is complete. 
\hfill  $\Box$ 

\vspace{\baselineskip}

 We also point out the validity of the following elementary but useful lemma.
 
\begin{lemma}\label{lem:coc-1ao}
 Let $n\in {\mathbb{N}}\setminus\{0\}$.  Let $\alpha\in]0,1]$. Let $\Omega$ be a bounded open  subset of ${\mathbb{R}}^{n}$. Let $X$ be a normed space. Let $L$ be a linear map from $C^{-1,\alpha}(\Omega)$ to $X$. Then $L$ is continuous if and only if the map
 \[
 L \circ  \Xi 
 \]
 is continuous on $C^{0,\alpha}(\overline{\Omega})^{n+1}$.
\end{lemma}
For a proof, we refer to \cite[Lem.~2.3]{La24c}.  We also mention the validity of the following approximation lemma.
\begin{lemma}\label{lem:aprc-1a}
 Let $\alpha\in]0,1]$. 
 Let $\Omega$ be a bounded open subset of ${\mathbb{R}}^n$ of class $C^{1,\alpha}$. 
If $f\in C^{-1,\alpha}(\overline{\Omega})$, then there exists a sequence $\{f_j\}_{j\in {\mathbb{N}}}$ in 
 $C^{\infty}(\overline{\Omega})$ such that
 \begin{equation}\label{lem:aprc-1a1}
 \sup_{j\in {\mathbb{N}}}\|f_j\|_{
C^{-1,\alpha}(\overline{\Omega}) 
 }<+\infty\,,\qquad
 \lim_{j\to\infty}f_j=f\quad\text{in}\ C^{-1,\beta}(\overline{\Omega})\quad \forall\beta\in]0,\alpha[\,.
 \end{equation}
\end{lemma}
{\bf Proof.} Let $(g_0,\dots,g_n )\in \left(C^{0,\alpha}(\overline{\Omega})\right)^{n+1}$ be such that
\begin{equation}
f=g_0+\sum_{s=1}^n\frac{\partial g_s}{\partial x_s}\,.
\end{equation}
A known approximation property implies that there exists a sequence $\{g_{s,j}\}_{j\in {\mathbb{N}}}$ in $C^{\infty}(\overline{\Omega})$   that converges to $g_s$ in the $ C^{0,\beta}(\overline{\Omega})$-norm for all $\beta\in]0,\alpha[$ and that is bounded in the $ C^{0,\alpha}(\overline{\Omega})$-norm, for each $s\in\{0,\dots,n\}$   (cf.~\cite[Lem.~A.3  of the Appendix]{La24c}). Since the map $\Xi$ from $\left(C^{0,\beta}(\overline{\Omega})\right)^{n+1}$ to $C^{-1,\beta}(\overline{\Omega})$ that takes a $(n+1)$-tuple $(\theta_0,\dots,\theta_n )$ to $\theta_0+\sum_{s=1}^n\frac{\partial \theta_s}{\partial x_s} $ is continuous and $\Xi$ is linear and continuous from $\left(C^{0,\alpha}(\overline{\Omega})\right)^{n+1}$ to $C^{-1,\alpha}(\overline{\Omega})$, we have
\begin{equation}\label{lem:aprc-1a3}
 \sup_{j\in {\mathbb{N}}}\|f_j\|_{
 C^{-1,\alpha}(\overline{\Omega}) 
 }<+\infty\,,\qquad
 \lim_{j\to\infty}f_j=f\quad\text{in}\ C^{-1,\beta}(\overline{\Omega})\quad \forall\beta\in]0,\alpha[\,,
 \end{equation}
where
\[
f_j\equiv g_{0,j}+\sum_{s=1}^n\frac{\partial g_{s,j}}{\partial x_s}\qquad\forall j\in {\mathbb{N}}\,.
\]
Hence, the proof is complete.\hfill  $\Box$ 

\vspace{\baselineskip}

We now define a linear functional ${\mathcal{I}}_{\Omega }$ on $C^{-1,\alpha}(\overline{\Omega})$ which extends the integration in $\Omega$ to all elements of   $C^{-1,\alpha}(\overline{\Omega})$ as in Dalla Riva, the author and Musolino \cite[Prop.~2.89]{DaLaMu21}.  
\begin{proposition}\label{prelim.wdtI}
  Let $\alpha\in]0,1]$. Let $\Omega$ be a bounded open Lipschitz subset of ${\mathbb{R}}^{n}$. Then there exists one and only one  linear and continuous  operator ${\mathcal{I}}_{\Omega}$ from the space $C^{-1,\alpha}(\overline{\Omega})$ to ${\mathbb{R}}$ such that
\begin{equation}\label{prelim.wdtI1}
{\mathcal{I}}_{\Omega}[f]=  \int_{\Omega}f_{0}\,dx+\int_{\partial\Omega}\sum_{j=1}^{n} (\nu_{\Omega})_{j}f_{j}\,d\sigma
\end{equation}
for all $f=  f_{0}+\sum_{j=1}^{n}\frac{\partial}{\partial x_{j}}f_{j}\in C^{-1,\alpha}(\overline{\Omega}) $. Moreover,
\[
{\mathcal{I}}_{\Omega}[f]=\int_{\Omega}f\,dx\qquad\forall f\in C^{0,\alpha}(\overline{\Omega})\,.
\]
\end{proposition}
We also exploit the following extension theorem, that enables to extend the elements of 
$C^{-1,\alpha}(\overline{\Omega})$, which are distributions in $\Omega$, to elements of the dual of $C^{1,\alpha}(\overline{\Omega})$. We do so by means of the following statement (see \cite[Prop.~3.1]{La24c} for a proof). 
\begin{proposition}\label{prop:nschext}
Let $\alpha\in]0,1]$. Let $\Omega$ be a bounded open Lipschitz subset of ${\mathbb{R}}^{n}$. Then there exists one and only one  linear and continuous extension operator $E^\sharp$ from $C^{-1,\alpha}(\overline{\Omega})$ to $\left(C^{1,\alpha}(\overline{\Omega})\right)'$ such that
 \begin{eqnarray}\label{prop:nschext2}
\lefteqn{
\langle E^\sharp[f],v\rangle 
}
\\ \nonumber
&&\ \
 =
\int_{\Omega}f_{0}v\,dx+\int_{\partial\Omega}\sum_{j=1}^{n} (\nu_{\Omega})_{j}f_{j}v\,d\sigma
 -\sum_{j=1}^{n}\int_{\Omega}f_{j}\frac{\partial v}{\partial x_j}\,dx
\quad \forall v\in C^{1,\alpha}(\overline{\Omega})
\end{eqnarray}
for all $f=  f_{0}+\sum_{j=1}^{n}\frac{\partial}{\partial x_{j}}f_{j}\in C^{-1,\alpha}(\overline{\Omega}) $. Moreover, 
\begin{equation}\label{prop:nschext1}
E^\sharp[f]_{|\Omega}=f\,, \ i.e.,\ 
\langle E^\sharp[f],v\rangle =\langle f,v\rangle \qquad\forall v\in {\mathcal{D}}(\Omega)
\end{equation}
for all $f\in C^{-1,\alpha}(\overline{\Omega})$ and
\begin{equation}\label{prop:nschext3}
\langle E^\sharp[f],v\rangle =\langle f,v\rangle \qquad\forall v\in C^{1,\alpha}(\overline{\Omega})
\end{equation}
for all $f\in C^{0,\alpha}(\overline{\Omega})$.
\end{proposition}

By  Proposition \ref{prop:nschext},   we know that the extension operator $E^\sharp$ defined as in (\ref{prop:nschext2}) satisfies condition (\ref{prop:nschext3}), but one may wonder whether such a choice can be considered as canonical. We answer by proving the following statement.
\begin{proposition}\label{prop:unschext}
 Let $\alpha\in]0,1]$. Let $\Omega$ be a bounded open Lipschitz subset of ${\mathbb{R}}^{n}$. 
 \begin{enumerate}
\item[(i)] If   $\beta\in]0,\alpha[$, then  $E^\sharp$  is continuous from $C^{-1,\alpha}(\overline{\Omega})$ with the norm of $C^{-1,\beta}(\overline{\Omega})$ to $\left(C^{1,\alpha}(\overline{\Omega})\right)'$ with the weak$^\ast$ topology.
\item[(ii)]
Let  $\tilde{E}^\sharp$ be a linear map from $C^{-1,\alpha}(\overline{\Omega})$ to $\left(C^{1,\alpha}(\overline{\Omega})\right)'$ that satisfies condition (\ref{prop:nschext3}) for all $f\in C^{\infty}(\overline{\Omega})$. If  there exists $\beta\in]0,\alpha[$ such that $\tilde{E}^\sharp$ 
 is continuous from $C^{-1,\alpha}(\overline{\Omega})$ with the norm of  $C^{-1,\beta}(\overline{\Omega})$ to  $\left(C^{1,\alpha}(\overline{\Omega})\right)'$ with the weak$^\ast$ topology, then $\tilde{E}^\sharp= E^\sharp$.
\end{enumerate}
\end{proposition}
{\bf Proof.} (i) By  Proposition \ref{prop:nschext}, there exists $E^\sharp_\beta\in {\mathcal{L}}\left(C^{-1,\beta}(\overline{\Omega}),\left(C^{1,\beta}(\overline{\Omega})\right)'\right)$ that satisfies conditions (\ref{prop:nschext2})--(\ref{prop:nschext3}) with $\beta$ instead of $\alpha$. Then $E^\sharp_\beta$ is continuous from $C^{-1,\beta}(\overline{\Omega})$    to $\left(C^{1,\beta}(\overline{\Omega})\right)'$ with the weak$^\ast$ topology. By equality  (\ref{prop:nschext2}), we have
\[
\langle E^\sharp[f],v\rangle=\langle E^\sharp_\beta[f],v\rangle\qquad\forall v\in C^{1,\alpha}(\overline{\Omega})
\subseteq C^{1,\beta}(\overline{\Omega})\,,
\]
for all $f\in C^{-1,\alpha}(\overline{\Omega})\left(\subseteq C^{-1,\beta}(\overline{\Omega})\right)$. Thus  if $v\in C^{1,\alpha}(\overline{\Omega})$, the map $\langle E^\sharp_\beta[\cdot],v\rangle$ is continuous from $C^{-1,\beta}(\overline{\Omega})$  to ${\mathbb{C}}$ and   $\langle E^\sharp [\cdot],v\rangle= \langle E^\sharp_\beta[\cdot]_{|C^{-1,\alpha}(\overline{\Omega})},v\rangle$  is continuous from $C^{-1,\alpha}(\overline{\Omega})$ with the norm of $C^{-1,\beta}(\overline{\Omega})$ to ${\mathbb{C}}$ and statement (i) holds true. 

(ii) Let $h\in  C^{-1,\alpha}(\overline{\Omega})$. By Lemma \ref{lem:aprc-1a}, there exists a sequence $\{h_l\}_{l\in {\mathbb{N}}}$   in $C^{\infty}(\overline{\Omega})$ as in (\ref{lem:aprc-1a1}). Then our continuity  assumption on  $\tilde{E}^\sharp$, the validity of condition (\ref{prop:nschext3})   for $\tilde{E}^\sharp$  with $f=h_l$ 
  and statement (i) imply that
\[
\langle  \tilde{E}^\sharp[h],v\rangle=\lim_{l\to\infty} \langle  \tilde{E}^\sharp[h_l],v\rangle=
\lim_{l\to\infty} \langle h_l,v\rangle=\lim_{l\to\infty} \langle E^\sharp [h_l],v\rangle=\langle E^\sharp[h],v\rangle
\]
for all $v\in C^{1,\alpha}(\overline{\Omega})$
and thus the proof is complete.\hfill  $\Box$ 

\vspace{\baselineskip}

In the specific case in which $\alpha=1$, we know that the elements of $C^{-1,1}(\overline{\Omega})$ are actually functions (cf.~Proposition \ref{prop:c-1ainfty}) and one can prove the following  simpler formula for the extension operator $E^\sharp$, that follows by applying the Divergence Theorem. 
\begin{proposition}\label{prop:nschexta=1}
 Let $\alpha\in]0,1]$. Let $\Omega$ be a bounded open Lipschitz subset of ${\mathbb{R}}^{n}$. If $f\in  C^{-1,1}(\overline{\Omega})$, then 
 \begin{equation}\label{prop:nschexta=12}
\langle  E^\sharp[f],v\rangle=\int_\Omega f v\,dx\qquad\forall v\in C^{1,1}(\overline{\Omega})\,,
\end{equation}
\textit{i.e.}, the extension operator $ E^\sharp$ from $C^{-1,1}(\overline{\Omega})$ to $\left(C^{1,1}(\overline{\Omega})\right)'$ coincides with  ${\mathcal{J}}$. (cf.~Lemma \ref{lem:cainclo} and Proposition \ref{prop:c-1ainfty}).
\end{proposition}
{\bf Proof.} By the membership of $f$ in $C^{-1,1}(\overline{\Omega})$, there exist $f_{j}\in C^{0,1}(\overline{\Omega})$ for all $j\in \{0,\dots,n\}$ such that
\[
f=f_{0}+\sum_{j=1}^{n}\frac{\partial}{\partial x_{j}}f_{j}\,.
\]
Then formula (\ref{prop:nschext2}) for $E^\sharp$ and  the Divergence Theorem (cf.~\textit{e.g.}, Mitrea, Mitrea and Mitrea \cite[Thm.~1.2.1]{MitMitMit22}), imply that
\begin{eqnarray*}
\lefteqn{
\langle E^\sharp[f],v\rangle 
}
\\ \nonumber
&&\ \
 =
\int_{\Omega}f_{0}v\,dx+\int_{\partial\Omega}\sum_{j=1}^{n} (\nu_{\Omega})_{j}f_{j}v\,d\sigma
 -\sum_{j=1}^{n}\int_{\Omega}f_{j}\frac{\partial v}{\partial x_j}\,dx
 \\ \nonumber
&&\ \
 =
\int_{\Omega}f_{0}v\,dx+\int_\Omega  \sum_{j=1}^n\frac{\partial }{\partial x_j}(f_{j}v)\,d\sigma
 -\sum_{j=1}^{n}\int_{\Omega}f_{j}\frac{\partial v}{\partial x_j}\,dx
  \\ \nonumber
&&\ \
 =\int_{\Omega}f_{0}v\,dx+\int_\Omega\sum_{j=1}^{n}\frac{\partial f_j }{\partial x_j}v\,dx
 =\int_{\Omega}fv\,dx\quad \forall v\in C^{1,1}(\overline{\Omega})\,,
 \end{eqnarray*}
 and thus the proof is complete.\hfill  $\Box$ 

\vspace{\baselineskip}

  \section{Preliminaries on the fundamental solution}
  \label{sec:prelfuso}

In order to analyze the volume potential, we need some more information on the fundamental solution $S_{ {\mathbf{a}} } $.  To do so, we introduce the fundamental solution $S_{n}$ of the Laplace operator. Namely, we set
\[
S_{n}(x)\equiv
\left\{
\begin{array}{lll}
\frac{1}{s_{n}}\ln  |x| \qquad &   \forall x\in 
{\mathbb{R}}^{n}\setminus\{0\},\quad & {\mathrm{if}}\ n=2\,,
\\
\frac{1}{(2-n)s_{n}}|x|^{2-n}\qquad &   \forall x\in 
{\mathbb{R}}^{n}\setminus\{0\},\quad & {\mathrm{if}}\ n>2\,,
\end{array}
\right.
\]
where $s_{n}$ denotes the $(n-1)$ dimensional measure of 
$\partial{\mathbb{B}}_{n}(0,1)$ and
we follow a formulation of Dalla Riva \cite[Thm.~5.2, 5.3]{Da13} and Dalla Riva, Morais and Musolino \cite[Thms.~3.1, 3.2, 5.5]{DaMoMu13}, that we state  as in paper  \cite[Cor.~4.2]{DoLa17} with Dondi (see also John~\cite{Jo55}, Miranda~\cite{Mi65} for homogeneous operators, and Mitrea and Mitrea~\cite[p.~203]{MitMit13}).   
\begin{proposition}
 \label{prop:ourfs} 
Let ${\mathbf{a}}$ be as in (\ref{introd0}), (\ref{ellip}), (\ref{symr}). 
Let $S_{ {\mathbf{a}} }$ be a fundamental solution of $P[{\mathbf{a}},D]$. 
Then there exist an invertible matrix $T\in M_{n}({\mathbb{R}})$ such that
\begin{equation}
\label{prop:ourfs0}
a^{(2)}=TT^{t}\,,
\end{equation}
 a real analytic function $A_{1}$ from $\partial{\mathbb{B}}_{n}(0,1)\times{\mathbb{R}}$ to ${\mathbb{C}}$ such that
 $A_{1}(\cdot,0)$ is odd,    $b_{0}\in {\mathbb{C}}$, a real analytic function $B_{1}$ from ${\mathbb{R}}^{n}$ to ${\mathbb{C}}$ such that $B_{1}(0)=0$, and a real analytic function $C $ from ${\mathbb{R}}^{n}$ to ${\mathbb{C}}$ such that
\begin{equation}
\label{prop:ourfs1}
S_{ {\mathbf{a}} }(x)
= 
\frac{1}{\sqrt{\det a^{(2)} }}S_{n}(T^{-1}x)
+|x|^{3-n}A_{1}(\frac{x}{|x|},|x|)
 +(B_{1}(x)+b_{0}(1-\delta_{2,n}))\ln  |x|+C(x)\,,
\end{equation}
for all $x\in {\mathbb{R}}^{n}\setminus\{0\}$,
 and such that both $b_{0}$ and $B_{1}$   equal zero
if $n$ is odd. Moreover, 
 \[
 \frac{1}{\sqrt{\det a^{(2)} }}S_{n}(T^{-1}x) 
 \]
is a fundamental solution for the principal part
  of $P[{\mathbf{a}},D]$.
\end{proposition}
In particular for the statement that $A_{1}(\cdot,0)$ is odd, we refer to
Dalla Riva, Morais and Musolino \cite[Thm.~5.5, (32)]{DaMoMu13}, where $A_{1}(\cdot,0)$ coincides with ${\mathbf{f}}_1({\mathbf{a}},\cdot)$ in that paper. Here we note that a function $A$ from $(\partial{\mathbb{B}}_{n}(0,1))\times{\mathbb{R}}$ to ${\mathbb{C}}$ is said to be real analytic provided that it has a real analytic extension   to an open neighbourhood of $(\partial{\mathbb{B}}_{n}(0,1))\times{\mathbb{R}}$ in 
${\mathbb{R}}^{n+1}$. Then we have the following elementary lemma    (cf.~\textit{e.g.},  \cite[Lem.~4.2]{La22d}).		 
\begin{lemma}\label{lem:anexsph}
 Let $n\in {\mathbb{N}}\setminus\{0,1\}$. A function $A$ from   $(\partial{\mathbb{B}}_{n}(0,1))\times{\mathbb{R}}$ to ${\mathbb{C}}$ is  real analytic if and only if the function $\tilde{A}$ from $({\mathbb{R}}^n\setminus\{0\}) \times{\mathbb{R}}$ defined by 
\begin{equation}\label{lem:anexsph1}
\tilde{A}(x,r)\equiv A(\frac{x}{|x|},r)\qquad\forall (x,r)\in ({\mathbb{R}}^n\setminus\{0\}) \times{\mathbb{R}}
\end{equation}
is real analytic.
 \end{lemma}
 
Then  one can prove the following formula for the Jacobian of the fundamental solution (see  Dondi and the author	 	 \cite[Lem.~4.3,  (4.8) and the following 2 lines]{DoLa17}). Here one should remember that $A_1(\cdot,0)$ is odd and that $b_0=0$ if $n$ is odd). 
\begin{proposition}
\label{prop:grafun}
 Let ${\mathbf{a}}$ be as in (\ref{introd0}), (\ref{ellip}), (\ref{symr}). Let $T\in M_{n}({\mathbb{R}})$  be as in (\ref{prop:ourfs0}). Let $S_{ {\mathbf{a}} }$ be a fundamental solution of $P[{\mathbf{a}},D]$. Let  $B_{1}$, $C$
 be as in Proposition \ref{prop:ourfs}. 
  Then there exists a real analytic function $A_{2}\equiv (A_{2,j})_{j=1,\dots,n}$ from $\partial{\mathbb{B}}_{n}(0,1)\times{\mathbb{R}}$ to ${\mathbb{C}}^{n}$ such that
\begin{eqnarray}
\label{prop:grafun1}
\lefteqn{
DS_{ {\mathbf{a}} }(x)=\frac{1}{ s_{n}\sqrt{\det a^{(2)} } }
|T^{-1}x|^{-n}x^{t}(a^{(2)})^{-1} 
}
\\ \nonumber
&&\qquad\qquad
+|x|^{2-n}A_{2}(\frac{x}{|x|},|x|)+DB_{1}(x)\ln |x|+DC(x)
\end{eqnarray}
for all $x\in {\mathbb{R}}^{n}\setminus\{0\}$. 
Moreover,   $A_2(\cdot,0)$ is even.
\end{proposition}
Next we introduce some notation. If $X$ and $Y$ are subsets of ${\mathbb{R}}^n$, then the symbol
\[
{\mathbb{D}}_{X\times Y}\equiv  \left\{
(x,y)\in X\times Y:\,x=y
\right\}
\] denotes the diagonal set of $X\times Y$  
 and we introduce the following class of   `potential type' kernels  (see also paper \cite{DoLa17} of the author and Dondi,  where such classes have been introduced in a form that generalizes those of Giraud \cite{Gi34}, Gegelia \cite{Ge67}, 
 Kupradze, Gegelia, Basheleishvili and 
 Burchuladze \cite[Chap.~IV]{KuGeBaBu79}).	 
\begin{definition}
Let $X$, $Y\subseteq {\mathbb{R}}^n$. Let $s_1$, $s_2$, $s_3\in {\mathbb{R}}$. We denote by the symbol ${\mathcal{K}}_{s_1, s_2, s_3} (X\times Y)$ the set of continuous functions $K$ from $(X\times Y)\setminus {\mathbb{D}}_{X\times Y}$ to ${\mathbb{C}}$ such that
 \begin{eqnarray*}
\lefteqn{
\|K\|_{  {\mathcal{K}}_{ s_1, s_2, s_3  }(X\times Y)  }
\equiv
\sup\biggl\{\biggr.
|x-y|^{ s_{1} }\vert K(x,y)\vert :\,(x,y)\in X\times Y, x\neq y
\biggl.\biggr\}
}
\\ \nonumber
&&\qquad\qquad\qquad
+\sup\biggl\{\biggr.
\frac{|x'-y|^{s_{2}}}{|x'-x''|^{s_{3}}}
\vert   K(x',y)- K(x'',y)  \vert :\,
\\ \nonumber
&&\qquad\qquad\qquad 
x',x''\in X, x'\neq x'', y\in Y\setminus {\mathbb{B}}_n(x',2|x'-x''|)
\biggl.\biggr\}<+\infty\,.
\end{eqnarray*}
\end{definition}
We now turn to compute the class of the (convolution) kernels that corresponds to the second order partial derivatives of the fundamental solution $S_{ {\mathbf{a}} }$ by means of the following statement.
\begin{proposition}
\label{prop:sedesa}
Let ${\mathbf{a}}$ be as in (\ref{introd0}), (\ref{ellip}), (\ref{symr}).  Let $S_{ {\mathbf{a}} }$ be a fundamental solution of $P[{\mathbf{a}},D]$. Let  $A_{2}\equiv (A_{2,j})_{j=1,\dots,n}$ , $B_{1}$, $C$ be as in Proposition \ref{prop:grafun} and formula
(\ref{prop:grafun1}). Let $G$ be a nonempty bounded subset of ${\mathbb{R}}^n$. 
Let 
\begin{equation}\label{prop:sedesa1}
k(x)=|x|^{2-n}A_{2}(\frac{x}{|x|},|x|)+DB_{1}(x)\ln |x|+DC(x)
\qquad\forall x\in {\mathbb{R}}^{n}\setminus\{0\}\,.
\end{equation}
Then the (convolution) kernel 
\[
\frac{\partial k}{\partial x_l}(x-y)\qquad\forall (x,y)\in (G\times G)\setminus {\mathbb{D}}_G
\]
belongs to $\left({\mathcal{K}}_{n-1,n,1}(G\times G)\right)^n$ for all $l\in\{1,\dots,n\}$. 
\end{proposition}
{\bf Proof.} Let $k_j$ denote the $j$-th component of $k$ for each $j\in\{1,\dots,n\}$. Let  $(\xi_1,\dots,\xi_n,r)$ denote  the
variable of $A_2$. 
Then we have
\begin{eqnarray}\label{prop:sedesa2}
\lefteqn{
\frac{\partial k_j}{\partial x_l}(x)=(2-n)|x|^{1-n}\frac{x_l}{|x|}A_{2,j}(\frac{x}{|x|},|x|)
}
\\ \nonumber
&&\qquad
+|x|^{2-n}\sum_{s=1}^n\frac{\partial A_{2,j}}{\partial\xi_s}(\frac{x}{|x|},|x|)
\left(
\frac{\delta_{sl}}{|x|}-\frac{x_sx_l}{|x|^3}
\right)
+|x|^{2-n}\frac{\partial A_{2,j}}{\partial r}(\frac{x}{|x|},|x|)\frac{x_l}{|x|}
\\ \nonumber
&&\qquad
+\frac{\partial^2 B_1}{\partial x_l\partial x_j}(x)\ln |x|+\frac{\partial B_1}{\partial x_j}(x)\frac{x_l}{|x|^2}+\frac{\partial^2 C}{\partial x_l\partial x_j}(x)
\qquad\forall x\in {\mathbb{R}}^{n}\setminus\{0\}\,,
\end{eqnarray}
for all $j$, $l\in\{1,\dots,n\}$.  Since $A_2$ is real analytic in $\partial {\mathbb{B}}_n(0,1)\times{\mathbb{R}}$,  Lemma 3.3 of Dondi and the author \cite{DoLa17} (see also Lemma 4.5 (i) of   \cite{La22d})  implies that the kernel
$A_{2,j}(\frac{x-y}{|x-y|},|x-y|)$ belongs to ${\mathcal{K}}_{0,1,1}(G\times G)$. Since the function $|\xi|^{1-n}\frac{\xi_l}{|\xi|}$ of the variable $\xi\in {\mathbb{R}}^n\setminus\{0\}$ is positively homogeneous of degree $-(n-1)$, Lemma 3.11 of \cite{La22d} implies that the kernel $|x-y|^{1-n}\frac{x_l-y_l}{|x-y|}$ is of class ${\mathcal{K}}_{n-1,n,1}(G\times G)$. Then the product Theorem 3.1
 (ii) of \cite{La22b} implies that the pointwise product is continuous from
\begin{equation}\label{prop:sedesa2a}
{\mathcal{K}}_{0,1,1}(G\times G)\times{\mathcal{K}}_{n-1,n,1}(G\times G)
\quad\text{to}\quad{\mathcal{K}}_{n-1,n,1}(G\times G) 
\end{equation}
and accordingly
\begin{equation}\label{prop:sedesa3}
(2-n)|x-y|^{1-n}\frac{x_l-y_l}{|x-y|}A_{2,j}(\frac{x-y}{|x-y|},|x-y|)
\in {\mathcal{K}}_{n-1,n,1}(G\times G)\,.
\end{equation}
We now consider the second addendum in the right hand side of equality (\ref{prop:sedesa2}). Since $\frac{\partial A_{2,j}}{\partial \xi_s}$ is real analytic in 
$\partial{\mathbb{B}}_n(0,1) \times {\mathbb{R}}$,    Lemma 3.3 of Dondi and the author \cite{DoLa17} (see also Lemma 4.5 (i) of   \cite{La22d})  implies that the kernel $\frac{\partial A_{2,j}}{\partial \xi_s}\left(\frac{x -y}{|x -y|},|x -y|\right)$ belongs to
${\mathcal{K}}_{0,1,1}(G\times G)$. Since the functions $|\xi|^{2-n}\left(
\frac{\delta_{sl}}{|\xi|}-\frac{\xi_s\xi_l}{|\xi|^3}
\right)$ of the variable $\xi\in {\mathbb{R}}^n\setminus\{0\}$ are positively homogeneous of degree $-(n-1)$, Lemma 3.11 of \cite{La22d} implies that  
\[
|x-y|^{2-n}\left(
\frac{\delta_{sl}}{|x-y|}-\frac{(x_s-y_s)(x_l-y_l) }{|x-y|^3}
\right)
\in{\mathcal{K}}_{n-1,n,1}(G\times G)\,.
\]
 Then the continuity of (\ref{prop:sedesa2a}) implies that
\begin{equation}\label{prop:sedesa4}
|x-y|^{2-n}\left(
\frac{\delta_{sl}}{|x-y|}-\frac{(x_s-y_s)(x_l-y_l) }{|x-y|^3}
\right)
\frac{\partial A_{2,j}}{\partial \xi_s}(\frac{x-y}{|x-y|},|x-y|)
\in {\mathcal{K}}_{n-1,n,1}(G\times G)\,.
\end{equation}
We now consider the third addendum in the right hand side of equality (\ref{prop:sedesa2}). Since $\frac{\partial A_{2,j}}{\partial r}$ is real analytic in 
$\partial{\mathbb{B}}_n(0,1) \times {\mathbb{R}}$,  Lemma 3.3 of Dondi and the author \cite{DoLa17} (see also Lemma 4.5 (i) of   \cite{La22d})  implies that
the kernel $\frac{\partial A_{2,j}}{\partial r}\left(\frac{x -y}{|x -y|},|x -y|\right)$ belongs to
${\mathcal{K}}_{0,1,1}(G\times G)$. Since the function    $|\xi|^{-(n-1)} \xi_l$  of the variable $\xi\in {\mathbb{R}}^n\setminus\{0\}$ is positively homogeneous of degree $-(n-2)$, Lemma 3.11 of \cite{La22d} implies that 
 \[
  |x-y|^{-(n-1)}(x_l-y_l)\in{\mathcal{K}}_{n-2,n-1,1}(G\times G)\,.
  \]
   Then the product Theorem 3.1
 (ii) of \cite{La22b} implies that the pointwise product is continuous from
\[
{\mathcal{K}}_{0,1,1}(G\times G)\times{\mathcal{K}}_{n-2,n-1,1}(G\times G)
\quad\text{to}\quad{\mathcal{K}}_{n-2,n-1,1}(G\times G)\,.
\]
Hence,
\begin{equation}\label{prop:sedesa5}
\frac{\partial A_{2,j}}{\partial r}\left(\frac{x -y}{|x -y|},|x -y|\right)|x-y|^{-(n-1)}(x_l-y_l)\in {\mathcal{K}}_{n-2,n-1,1}(G\times G)\,.
\end{equation}
We now consider the fourth addendum in the right hand side of equality (\ref{prop:sedesa2}). 
Since $B_1$ is analytic, Lemma 4.5 (ii) of \cite{La22d} implies that the kernel $\frac{\partial^2B_1}{\partial x_l\partial x_j}(x-y)$ belongs to
${\mathcal{K}}_{0,0,1}(G\times G)$ that is contained in 
${\mathcal{K}}_{0,1,1}(G\times G)$ (cf.~Proposition 3.2 (ii)  of \cite{La22b}). By Lemma 4.5 (iii) of \cite{La22d}  and by the embedding  Proposition 3.2 (ii) of \cite{La22b}, we have
\[
\ln |x-y|\in {\mathcal{K}}_{\epsilon,1,1}(G\times G)\subseteq{\mathcal{K}}_{\epsilon,1+\epsilon,1}(G\times G)\qquad\forall\epsilon\in]0,1[\,.
\]
Then the product Theorem 3.1
 (ii) of \cite{La22b}  implies that the pointwise product is continuous from
\[
{\mathcal{K}}_{0,1,1}(G\times G)\times{\mathcal{K}}_{\epsilon,1+\epsilon,1}(G\times G)
\quad\text{to}\quad{\mathcal{K}}_{\epsilon,1+\epsilon,1}(G\times G)\qquad\forall\epsilon\in]0,1[ 
\]
and accordingly
\begin{equation}\label{prop:sedesa6}
\frac{\partial^2B_1}{\partial x_l\partial x_j}(x-y)\ln |x-y|\in {\mathcal{K}}_{\epsilon,1+\epsilon,1}(G\times G)\qquad\forall\epsilon\in]0,1[\,.
\end{equation}
We now consider the fifth addendum in the right hand side of equality (\ref{prop:sedesa2}). 
Since $B_1$ is analytic, Lemma 4.5 (ii) of \cite{La22d} implies that the kernel $\frac{\partial B_1}{\partial x_j}(x-y)$ belongs to ${\mathcal{K}}_{0,0,1}(G\times G)$ that is contained in ${\mathcal{K}}_{0,1,1}(G\times G)$ (cf.~Proposition 3.5 (ii)  of \cite{La22b}). 
Since the function    $|\xi|^{-2} \xi_l$  of the variable $\xi\in {\mathbb{R}}^n\setminus\{0\}$ is positively homogeneous of degree $-1$, Lemma 3.11 of \cite{La22d} implies that 
the kernels $|x-y|^{-2}(x_l-y_l)$  are of class ${\mathcal{K}}_{1,2,1}(G\times G)$. Then the product Theorem 3.1
 (ii) of \cite{La22b} implies that the pointwise product is continuous from
\[
{\mathcal{K}}_{0,1,1}(G\times G)\times{\mathcal{K}}_{1,2,1}(G\times G)
\quad\text{to}\quad{\mathcal{K}}_{1,2,1}(G\times G) 
\]
and accordingly
\begin{equation}\label{prop:sedesa7}
\frac{\partial B_1}{\partial x_j}(x-y)|x-y|^{-2}(x_l-y_l)\in {\mathcal{K}}_{1,2,1}(G\times G)\,.
\end{equation}
We now consider the sixth addendum in the right hand side of equality (\ref{prop:sedesa2}). Since $C$ is analytic, Lemma 4.5 (ii) of \cite{La22d} implies that the kernel
$\frac{\partial^2C}{\partial x_l\partial x_j}(x-y)$ belongs to
${\mathcal{K}}_{0,0,1}(G\times G)$ that is contained in ${\mathcal{K}}_{0,1,1}(G\times G)$ (cf.~Proposition 3.2 (ii)  of \cite{La22b}). Then
\begin{equation}\label{prop:sedesa8}
\frac{\partial^2C}{\partial x_l\partial x_j}(x-y)\in {\mathcal{K}}_{0,1,1}(G\times G)\,.
\end{equation}
Thus we have proved that each addendum in the right hand side of equality (\ref{prop:sedesa2}) is contained in one of the following classes
\begin{eqnarray*}
&&{\mathcal{K}}_{n-1,n,1}(G\times G)\,,\qquad {\mathcal{K}}_{n-2,n-1,1}(G\times G)\,,
\\
&&{\mathcal{K}}_{\epsilon,1+\epsilon,1}(G\times G)\qquad\forall\epsilon\in]0,1[\,,\qquad {\mathcal{K}}_{1,2,1}(G\times G)\,,\qquad
{\mathcal{K}}_{0,1,1}(G\times G)
\,.
\end{eqnarray*}
Now the embedding Proposition 3.2 of \cite{La22b} implies that each of such classes is contained in ${\mathcal{K}}_{n-1,n,1}(G\times G)$ and thus the proof is complete.\hfill  $\Box$ 

\vspace{\baselineskip}

\section{The distributional  volume potential}
  \label{sec:dvopo}

\begin{definition}\label{defn:dvpsa}
 Let $\alpha\in]0,1]$. 
 Let   $\Omega$ be a bounded open subset of ${\mathbb{R}}^{n}$. 
  Let ${\mathbf{a}}$ be as in (\ref{introd0}), (\ref{ellip}), (\ref{symr}).   Let $S_{ {\mathbf{a}} }$ be a fundamental solution of $P[{\mathbf{a}},D]$.
  If $\mu\in (C^{1,\alpha}(\overline{\Omega}))'$, then the (distributional) volume potential relative to $S_{ {\mathbf{a}} }$ and $\mu$ is the distribution
\[
{\mathcal{P}}_\Omega[\mu]=(r_{|\overline{\Omega}}^t\mu)\ast S_{ {\mathbf{a}} }\in {\mathcal{D}}'({\mathbb{R}}^n)\,.
\]
\end{definition}
By the definition of convolution of distributions, we have
\begin{eqnarray*}
\lefteqn{
\langle (r_{|\overline{\Omega}}^t\mu)\ast S_{ {\mathbf{a}} },\varphi\rangle =\langle r_{|\overline{\Omega}}^t\mu(y),\langle S_{ {\mathbf{a}} }(\eta),\varphi(y+\eta)\rangle \rangle 
}
\\ \nonumber
&&\qquad
=\langle r_{|\overline{\Omega}}^t\mu(y),\int_{{\mathbb{R}}^{n}}S_{ {\mathbf{a}} }(\eta)\varphi(y+\eta)\,d\eta\rangle 
=\langle r_{|\overline{\Omega}}^t\mu(y),\int_{{\mathbb{R}}^{n}}S_{ {\mathbf{a}} }(x-y)\varphi(x)\,dx\rangle 
\end{eqnarray*}
for all $\varphi\in{\mathcal{D}}({\mathbb{R}}^n)$. In general, $(r_{|\overline{\Omega}}^t\mu)\ast S_{ {\mathbf{a}} }$ is not a function, \textit{i.e.} $(r_{|\overline{\Omega}}^t\mu)\ast S_{ {\mathbf{a}} }$ is not a distribution that is associated to a locally integrable function in ${\mathbb{R}}^n$. 
However, this is the case if for example  $\mu$ is associated to a function of $ L^\infty(\Omega)$, \textit{i.e.}, $\mu={\mathcal{J}}[f]$ with $f\in  L^\infty(\Omega)$
(see Lemma \ref{lem:cainclo} with any choice of  $\alpha\in]0,1]$). Indeed, 
\begin{eqnarray*}
\lefteqn{
\langle (r_{|\overline{\Omega}}^t\mu)\ast S_{ {\mathbf{a}} },\varphi\rangle  =\langle (r_{|\overline{\Omega}}^t{\mathcal{J}}[f])\ast S_{ {\mathbf{a}} },\varphi\rangle  
}
\\ \nonumber
&&\qquad
=\langle r_{|\overline{\Omega}}^t{\mathcal{J}}[f](y),\int_{\Omega}S_{ {\mathbf{a}} }(x-y)\varphi(x)\,dx\rangle 
\\ \nonumber
&&\qquad
=\langle {\mathcal{J}}[f](y),r_{|\overline{\Omega}}\int_{{\mathbb{R}}^{n}}S_{ {\mathbf{a}} }(x-y)\varphi(x)\,dx\rangle 
\\ \nonumber
&&\qquad
=\int_{\Omega}f(y)\int_{{\mathbb{R}}^{n}}S_{ {\mathbf{a}} }(x-y)\varphi(x)\,dx\,dy
=\int_{{\mathbb{R}}^{n}}\int_{\Omega}S_{ {\mathbf{a}} }(x-y)f(y)\,dy\varphi(x)\,dx
\\ \nonumber
&&\qquad
=\langle \int_{\Omega}S_{ {\mathbf{a}} }(x-y)f(y)\,dy,\varphi(x)\rangle 
\end{eqnarray*}
for all $\varphi\in{\mathcal{D}}({\mathbb{R}}^n)$ and thus 
 the (distributional) volume potential relative to $S_{ {\mathbf{a}} }$ and $\mu$ is associated to the function
\begin{equation}\label{prop:dvpsa1}
\int_{\Omega}S_{ {\mathbf{a}} }(x-y)f(y)\,dy\qquad{\mathrm{a.a.}}\ x\in {\mathbb{R}}^n\,,
\end{equation}
that is locally integrable in ${\mathbb{R}}^n$  and that with some abuse of notation we still denote by the symbol ${\mathcal{P}}_\Omega[S_{ {\mathbf{a}} },{\mathcal{J}}[f]]$  or even more  simply by the symbol  ${\mathcal{P}}_\Omega[S_{ {\mathbf{a}} }, f]$.  We also note that  under the assumptions of Definition \ref{defn:dvpsa}, classical   properties of the convolution of distributions imply that
 \begin{equation}\label{prop:dvpsa2}
P[{\mathbf{a}},D] \left[(r_{|\overline{\Omega}}^t\mu)\ast S_{ {\mathbf{a}} }\right]
 =(r_{|\overline{\Omega}}^t\mu)\ast (P[{\mathbf{a}},D][S_{ {\mathbf{a}} }])=
 (r_{|\overline{\Omega}}^t\mu)\ast\delta_0=(r_{|\overline{\Omega}}^t\mu) \quad\text{in}\ {\mathcal{D}}'({\mathbb{R}}^n)\,,
 \end{equation} 
 where $\delta_0$ is the Dirac measure with mass at $0$.  We now present a classical formula for the function that represents the restriction of  the distributional volume  potential $(r_{|\overline{\Omega}}^t\mu)\ast S_{ {\mathbf{a}} }$ to ${\mathbb{R}}^n\setminus {\mathrm{supp}}\,(r_{|\overline{\Omega}}^t\mu)$ (and thus to 
 ${\mathbb{R}}^n\setminus \overline{\Omega}$) by means of the following statement. For the convenience of the reader, we include a proof.
\begin{proposition}\label{prop:dvpfun} 
 Let ${\mathbf{a}}$ be as in (\ref{introd0}), (\ref{ellip}), (\ref{symr}).   Let $S_{ {\mathbf{a}} }$ be a fundamental solution of $P[{\mathbf{a}},D]$. Let     $\tau\in {\mathcal{D}}'({\mathbb{R}}^n)$ be a distribution with compact support  ${\mathrm{supp}}\,\tau$.  Then the real valued function $\theta$ from ${\mathbb{R}}^n\setminus {\mathrm{supp}}\,\tau$ that is defined by 
\begin{equation}\label{prop:dvpfun1}
\theta(x)\equiv\langle \tau (y),S_{ {\mathbf{a}} }(x-y)\rangle \qquad\forall x\in {\mathbb{R}}^n\setminus {\mathrm{supp}}\,\tau
\end{equation}
is of class $C^\infty$ and the restriction of   $\tau\ast S_{ {\mathbf{a}} }$ to ${\mathbb{R}}^n\setminus {\mathrm{supp}}\,\tau$ is associated to the function $\theta$. Namely,
\begin{equation}\label{prop:dvpfun2}
\langle \tau\ast S_{ {\mathbf{a}} },\varphi\rangle =\int_{{\mathbb{R}}^n\setminus{\mathrm{supp}}\,\tau}
 \langle \tau(y),S_{ {\mathbf{a}} }(x-y)\rangle \varphi(x)\,dx
 \quad\forall \varphi\in {\mathcal{D}}({\mathbb{R}}^n\setminus {\mathrm{supp}}\,\tau)\,.
 \end{equation}
  [Here we note that  the symbol 
 $\langle \tau (y),S_{ {\mathbf{a}} }(x-y)\rangle $ in (\ref{prop:dvpfun1}) means  
 \[
 \langle \tau (y),\omega(y)S_{ {\mathbf{a}} }(x-y)\rangle \,,
 \]
  where 
 $\omega\in {\mathcal{D}}({\mathbb{R}}^n\setminus \{x\})$ and $\omega$ equals $1$ in an open neighborhood of  ${\mathrm{supp}}\,\tau$.] Moreover, $P[{\mathbf{a}},D][\theta]=0$ in ${\mathbb{R}}^n\setminus {\mathrm{supp}}\,\tau$.
 \end{proposition}
{\bf Proof.}   Since $\tau$ is a distribution in ${\mathbb{R}}^n$ with compact support and $S_{ {\mathbf{a}} }(x-\cdot)$ is of class $C^\infty$ in ${\mathbb{R}}^n\setminus\{x\}$ for all $x\in {\mathbb{R}}^n\setminus {\mathrm{supp}}\,\tau$,  the differentiablity theorem for distributions with compact support in ${\mathbb{R}}^n$
 applied to test functions depending on a parameter implies that  the function $\theta$ is of class $C^\infty$ in ${\mathbb{R}}^n\setminus {\mathrm{supp}}\,\tau$ (cf.~\textit{e.g.}, Treves \cite[Thm.~27.2]{Tr67}). We now fix $\varphi\in {\mathcal{D}}({\mathbb{R}}^n\setminus {\mathrm{supp}}\,\tau)$ and we prove equality (\ref{prop:dvpfun2}). 
 
 Let $\Omega^\sharp$ be an open neighborhood of $ {\mathrm{supp}}\,\tau$ such that $\overline{\Omega^\sharp}\cap {\mathrm{supp}}\,\varphi=\emptyset$. By the known sequential density of ${\mathcal{D}}(\Omega^\sharp)$ in the space of compactly supported distributions in $\Omega^\sharp$, there exists a sequence $\{\tau_j\}_{j\in {\mathbb{N}} }$ in ${\mathcal{D}}(\Omega^\sharp)$ such that
\begin{equation}\label{prop:dvpfun3}
\lim_{j\to\infty}\tau_j=\tau\qquad\text{in} \ (C^\infty(\Omega^\sharp))'_b\,,
\end{equation}
and accordingly in $(C^\infty({\mathbb{R}}^n))'_b$, where $(C^\infty(\Omega^\sharp))'_b$ and $(C^\infty({\mathbb{R}}^n))'_b$ denote  the dual of $C^\infty(\Omega^\sharp)$ with the topology of uniform convergence on the bounded subsets of $C^\infty(\Omega^\sharp)$
and the dual of $C^\infty({\mathbb{R}}^n)$ with the topology of uniform convergence on the bounded subsets of $C^\infty({\mathbb{R}}^n)$, respectively (cf.~\textit{e.g.}, Treves \cite[Thm.~28.2]{Tr67}. See Treves \cite[Chapt.~10, Ex.~I, Chapt.~14]{Tr67} for the definition of  topology of $C^\infty(\Omega^\sharp)$  and of bounded subsets of $C^\infty(\Omega^\sharp)$). 

Then the above mentioned   differentiablity theorem for distributions with compact support in ${\mathbb{R}}^n$
 applied to test functions depending on a parameter implies that  the function
 $\langle \tau_j(y),S_{ {\mathbf{a}} }(\cdot-y)\rangle $ is of class $C^\infty$ in ${\mathbb{R}}^n\setminus \overline{\Omega^\sharp}$ for each $j\in{\mathbb{N}}$. 
 By  the definition of convolution and the convergence of (\ref{prop:dvpfun3}) in $(C^\infty({\mathbb{R}}^n))'_b$ and accordingly in $(C^\infty(\Omega^\sharp))'_b$, we have
 \begin{eqnarray}\label{prop:dvpfun4}
\lefteqn{
\langle \tau\ast S_{ {\mathbf{a}} },\varphi\rangle 
=\langle \tau(y),\langle S_{ {\mathbf{a}} }(\eta),\varphi(y+\eta)\rangle \rangle 
}
\\ \nonumber
&&\qquad\qquad
=\lim_{j\to\infty}\langle \tau_j(y),\langle S_{ {\mathbf{a}} }(\eta),\varphi(y+\eta)\rangle \rangle 
\\ \nonumber
&&\qquad\qquad
=\lim_{j\to\infty}\int_{{\mathbb{R}}^n}\tau_j(y)\int_{{\mathbb{R}}^n}S_{ {\mathbf{a}} }(\eta) \varphi(y+\eta)\,d\eta\,dy
\\ \nonumber
&&\qquad\qquad
=\lim_{j\to\infty}\int_{{\mathbb{R}}^n}\tau_j(y)\int_{{\mathbb{R}}^n}S_{ {\mathbf{a}} }(x-y) \varphi(x)\,dx\,dy 
\\ \nonumber
&&\qquad\qquad
=\lim_{j\to\infty}\int_{{\mathbb{R}}^n} \int_{{\mathbb{R}}^n}\tau_j(y)S_{ {\mathbf{a}} }(x-y)\,dy  \varphi(x)\,dx 
\\ \nonumber
&&\qquad\qquad
=\lim_{j\to\infty}\int_{{\mathbb{R}}^n} \langle \tau_j(y),S_{ {\mathbf{a}} }(x-y)\rangle  \varphi(x)\,dx \,.
\end{eqnarray}
Next we turn to show that the sequence $\{\langle \tau_j(y),S_{ {\mathbf{a}} }(x-y)\rangle \}_{j\in{\mathbb{N}} }$ converges uniformly to $\langle \tau(y),S_{ {\mathbf{a}} }(x-y)\rangle $ in $x\in{\mathrm{supp}}\,\varphi$. 
   Since $\Omega^\sharp$ has a strictly  positive distance from $ {\mathrm{supp}}\,\varphi$, the  set
\[
\{S_{ {\mathbf{a}} }(x-\cdot):\,x\in {\mathrm{supp}}\,\varphi\}
\]
 is bounded in $C^\infty(\Omega^\sharp)$ and accordingly
 \[
 \lim_{j\to\infty}\langle \tau_j(y),S_{ {\mathbf{a}} }(x-y)\rangle =\langle \tau,S_{ {\mathbf{a}} }(x-y)\rangle 
 \]
 uniformly in $x\in {\mathrm{supp}}\,\varphi$.  Hence,
 \[
 \lim_{j\to\infty}\int_{{\mathbb{R}}^n}\langle \tau_j(y),S_{ {\mathbf{a}} }(x-y)\rangle \varphi(x)\,dx=
 \int_{{\mathbb{R}}^n}\langle \tau(y),S_{ {\mathbf{a}} }(x-y)\rangle \varphi(x)\,dx
 \]
   and equality (\ref{prop:dvpfun4}) implies that equality (\ref{prop:dvpfun2}) holds true. Moreover, known properties of the convolution imply that
   \begin{equation}\label{prop:dvpsa2a}
 P[{\mathbf{a}},D]\left[\tau\ast S_{ {\mathbf{a}} }\right]
 =\tau\ast (P[{\mathbf{a}},D][S_{ {\mathbf{a}} }])=
\tau\ast\delta_0=\tau \quad\text{in}\ {\mathcal{D}}'({\mathbb{R}}^n)\,.
 \end{equation} 
 Since $\tau$ vanishes in ${\mathbb{R}}^n\setminus {\mathrm{supp}}\,\tau$ and $P[{\mathbf{a}},D]$ is elliptic,  the function $\theta$ that represents the  restriction of $\tau\ast S_{ {\mathbf{a}} }$ to ${\mathbb{R}}^n\setminus {\mathrm{supp}}\,\tau$  is real analytic and $P[{\mathbf{a}},D][\theta]=0$  in ${\mathbb{R}}^n\setminus {\mathrm{supp}}\,\tau$.\hfill  $\Box$ 

\vspace{\baselineskip}

\section{A technical lemma on an integral operator}
\label{sec:teleinop}

We first introduce two (known) normed spaces of positively homogeneous functions.
If $n\in {\mathbb{N}}\setminus\{0\}$, $m\in {\mathbb{N}}$, $h\in {\mathbb{R}}$, $\alpha\in ]0,1]$, then we set 
\begin{equation}\label{volume.klhomog}
{\mathcal{K}}^{m,\alpha}_h \equiv\biggl\{
k\in C^{m,\alpha}_{ {\mathrm{loc}}}({\mathbb{R}}^n\setminus\{0\}):\, k\ {\text{is\ positively\ homogeneous\ of \ degree}}\ h
\biggr\}\,,
\end{equation}
where $C^{m,\alpha}_{ {\mathrm{loc}}}({\mathbb{R}}^n\setminus\{0\})$ denotes the set of functions of 
$C^{m}({\mathbb{R}}^n\setminus\{0\})$ whose restriction to $\overline{\Omega}$ is of class $C^{m,\alpha}(\overline{\Omega})$ for all bounded open subsets $\Omega$ of ${\mathbb{R}}^n$ such that
$\overline{\Omega}\subseteq {\mathbb{R}}^n\setminus\{0\}$
and we set
\[
\|k\|_{ {\mathcal{K}}^{m,\alpha}_h}\equiv \|k\|_{C^{m,\alpha}(\partial{\mathbb{B}}_n(0,1))}\qquad\forall k\in {\mathcal{K}}^{m,\alpha}_h\,.
\]
We can easily verify that $ \left({\mathcal{K}}^{m,\alpha}_h , \|\cdot\|_{ {\mathcal{K}}^{m,\alpha}_h}\right)$ is a Banach space and we  consider the closed subspaces
\begin{eqnarray}\label{volume.k01e0}
{\mathcal{K}}^{m,\alpha}_{h;o}&\equiv& \left\{
k\in {\mathcal{K}}^{m,\alpha}_h:\,k\ {\mathrm{is\ odd}} 
\right\}\,,
\\ \nonumber
{\mathcal{K}}^{m,\alpha}_{h;e,0}&\equiv& \left\{
k\in {\mathcal{K}}^{m,\alpha}_h:\,k\ {\mathrm{is\ even}}, \int_{\partial{\mathbb{B}}_n(0,1)}k\,d\sigma=0
\right\}
\end{eqnarray}
of ${\mathcal{K}}^{m,\alpha}_h$. Next we introduce the following known lemma on the maximal function associated to a convolution kernel in the specific case in which 
$k$ is even and has   integral equal to zero on the unit sphere and 
$\Omega$ is of class $C^{1,\alpha}$ (\textit{cf.} Majda and Bertozzi \cite[Prop. 8.12, pp.~348--350]{MaBe02}). For a proof, we refer to  the  proof due to  Mateu, Orobitg, and Verdera \cite[estimate of $(IV)_\delta$, p.~408]{MaOrVe08} (see also \cite[Lem.~7.11]{DaLaMu21}).
\begin{lemma}\label{lem:maorve} 
Let   $\alpha\in]0,1]$. 
 Let $\Omega$ be a bounded open subset of  ${\mathbb{R}}^n$ of class $C^{1,\alpha}$. Then there exists $c^*_\Omega\in]0,+\infty[$ such that\index{$c^*_\Omega$}
\begin{equation}\label{lem:maorve1}
\sup_{x\in {\mathbb{R}}^n}
\sup_{\rho\in ]0,+\infty[}
\left|
\int_{\Omega\setminus{\mathbb{B}}_{n}(x,\rho)}
k(x-y)\,dy
\right|
\leq c^*_\Omega\left\| k \right\|_{{\mathcal{K}}^{0,1}_{-n }}
\quad\forall k\in {\mathcal{K}}^{0,1}_{-n; e,0}\,.
\end{equation}
\end{lemma}
It is also known that if $k\in {\mathcal{K}}^{1,1}_{-(n-1);o}$, then its first order partial derivatives belong to ${\mathcal{K}}^{0,1}_{-n; e,0}$. Namely, the following holds. For a proof we refer for example to \cite[Lemmas 4.13, 7.12]{DaLaMu21}.
\begin{lemma}\label{lem:codepoh}
 Let $n\in{\mathbb{N}}\setminus\{0\}$, $j\in\{1,\dots,n\}$. Then the linear operator from ${\mathcal{K}}^{1,1}_{-(n-1);o}$ to ${\mathcal{K}}^{0,1}_{-n; e,0}$ that takes $k$ to $\frac{\partial k}{\partial x_j}$ is continuous.
\end{lemma}
Next we deduce the following extension of a statement of Dalla Riva, the author and Musolino  \cite[Prop.~7.15]{DaLaMu21} 
 by means of an abstract result of \cite[Prop.~6.3 (ii), (b)]{La22a} and by \cite[Lem.~3.11]{La22d}.
\begin{proposition}\label{prop:gjhc}
 Let   $\alpha\in]0,1]$,  $l\in\{1,\dots,n\}$. Let $\Omega$ be a bounded open subset of ${\mathbb{R}}^{n}$ of class $C^{1,\alpha}$. Let $r\in]0,+\infty[$ be such that $\overline{\Omega}\subseteq {\mathbb{B}}_n(0,r)$. 
 Let  
 \begin{equation}\label{prop:gjhc1}
G_l[k,\psi](x)\equiv\int_\Omega\frac{\partial k}{\partial x_l}(x-y)(\psi(y)-\psi(x))\,dy\qquad\forall x\in \overline{{\mathbb{B}}_n(0,r)} 
\end{equation}
for all $(k,\psi)\in {\mathcal{K}}^{1,1}_{-(n-1);o}\times C^{0,\alpha}(\overline{{\mathbb{B}}_n(0,r)})$. Then the following statements hold. 
\begin{enumerate}
\item[(i)] If $\alpha\in]0,1[$, then the bilinear map $G_l[\cdot,\cdot]$ from ${\mathcal{K}}^{1,1}_{-(n-1);o}\times C^{0,\alpha}(\overline{{\mathbb{B}}_n(0,r)})$ to $
 C^{0,\alpha}(\overline{{\mathbb{B}}_n(0,r)})$ that is delivered by the formula (\ref{prop:gjhc1}) is continuous.
\item[(ii)] If $\alpha=1$, then the bilinear map $G_l[\cdot,\cdot]$ from ${\mathcal{K}}^{1,1}_{-(n-1);o}\times C^{0,\alpha}(\overline{{\mathbb{B}}_n(0,r)})$ to $
 C^{0,\omega_1(\cdot)}(\overline{{\mathbb{B}}_n(0,r)})$ that is delivered by the formula (\ref{prop:gjhc1}) is continuous.
\end{enumerate}
\end{proposition}
{\bf Proof.} We first set
\[
X\equiv \overline{{\mathbb{B}}_n(0,r)}  \,,\qquad Y\equiv \Omega\,.
\]
Then we obviously have
\begin{eqnarray*}
\lefteqn{m_n(
({\mathbb{B}}_n(x,\rho_2)\setminus {\mathbb{B}}_n(x,\rho_1))\cap \Omega
)\leq m_n({\mathbb{B}}_n(0,1))(\rho_2^n-\rho_1^n)
}
\\ \nonumber
&&\qquad\qquad\qquad\qquad\qquad\qquad\qquad\qquad\forall x\in X,  \rho_1, \rho_2\in [0,+\infty[\ \text{with}\  \rho_1<\rho_2
\end{eqnarray*}
and accordingly  $Y$ is strongly upper $n$-Ahlfors regular with respect to $X$  in the sense of  \cite[(1.5)]{La22a}. We plan to apply  an abstract result of \cite[Prop.~6.3 (ii) (b), (bb)]{La22a}. Thus we note that if we set
\[
\upsilon_Y\equiv n\,,\quad s_1\equiv n\,,\quad s_2\equiv n+1\,,\quad s_3\equiv 1\,,
\]
then we have
\[
\upsilon_Y\in]0,+\infty[\,, \quad
s_1\in[\alpha,\upsilon_Y+\alpha[\,,\quad
s_2\in [\alpha,+\infty[\,,\quad s_3\in]0,1]
\qquad\text{if}\ \alpha\in]0,1]\,,
\]
\begin{eqnarray*}
&&s_2-\alpha=n+1-\alpha>n=\upsilon_Y\,,\quad s_2=n+1<n+\alpha+1=\upsilon_Y+\alpha+s_3\,,
\\
&&C^{0,\min\{\alpha,\upsilon_Y+s_3+\alpha-s_2\}}_b(X) = C^{0,\alpha}_b(X)\qquad\text{if}\ \alpha\in]0,1[\,,
\end{eqnarray*}
and
\[
s_2-\alpha=n=\upsilon_Y
\,,\qquad C^{0,\max\{r^\alpha,\omega_1(r)\}}_b(X)
 = C^{0, \omega_1(r) }_b(X)
\qquad\text{if}\ \alpha=1
\]
Then  \cite[Prop.~6.3 (ii) (b), (bb), Defn.~6.2]{La22a} and \cite[Lem.~3.11]{La22d} imply that there exist $c_{\alpha,0}$, $c_{\alpha,1}\in ]0,+\infty[$ such that
\begin{eqnarray}\label{prop:gjhc2}
\lefteqn{
\left.
\begin{array}{ll}
\text{if}\ \alpha\in]0,1[
 	&  \|G_l[k,\psi]\|_{
C^{0,\alpha} ( \overline{{\mathbb{B}}_n(0,r)} )}
\\
\text{if}\ \alpha=1 & \|G_l[k,\psi]\|_{
C^{0,\omega_1(\cdot)} ( \overline{{\mathbb{B}}_n(0,r)} )
}   
\end{array}
\right\} }
\\ \nonumber
&& 
\leq c_{\alpha,0}\|\psi\|_{C^{0,\alpha}_b(X\cup Y)}
 \biggl(
\left\|\frac{\partial k}{\partial x_l}(x-y)\right\|_{
{\mathcal{K}}_{n,n+1,1}(X\times Y)
}
\\ \nonumber
&&\qquad 
+\sup_{x\in \overline{{\mathbb{B}}_n(0,r)}}\sup_{\rho\in]0,+\infty[}
\left|
\int_{\Omega\setminus {\mathbb{B}}_n(x,\rho)}
\frac{\partial k}{\partial x_l}(x-y)\,dy
\right|\biggr) 
\\ \nonumber
&& \leq c_{\alpha,1}\|\psi\|_{C^{0,\alpha}(\overline{{\mathbb{B}}_n(0,r)})}
\\ \nonumber
&&\qquad
\times\left(
\left\|\frac{\partial k}{\partial x_l}\right\|_{
{\mathcal{K}}^{0,1}_{-n}
}
+\sup_{x\in \overline{{\mathbb{B}}_n(0,r)}	}\sup_{\rho\in]0,+\infty[}
\left|
\int_{\Omega\setminus {\mathbb{B}}_n(x,\rho)}
\frac{\partial k}{\partial x_l}(x-y)\,dy
\right|\right) 
\end{eqnarray} 
for all $(k,\psi)\in {\mathcal{K}}^{1,1}_{-(n-1);o}\times C^{0,\alpha}(\overline{{\mathbb{B}}_n(0,r)})$.  Then inequality (\ref{prop:gjhc2}) and Lemmas \ref{lem:maorve}, \ref{lem:codepoh} imply that $G_l$ is bilinear and continuous from 
\[
{\mathcal{K}}^{1,1}_{-(n-1);o}\times C^{0,\alpha}(\overline{{\mathbb{B}}_n(0,r)})\quad\text{to}\quad 
\left\{
\begin{array}{ll}
C^{0,\alpha} ( \overline{{\mathbb{B}}_n(0,r)}) 
&\text{if}\ \alpha\in]0,1[\,,
\\
C^{0,\omega_1(\cdot)} ( \overline{{\mathbb{B}}_n(0,r)}) 
&\text{if}\ \alpha=1
\end{array}
\right.
\]
and thus the proof is complete.\hfill  $\Box$ 

\vspace{\baselineskip}

We are now ready to prove  the following extension of a statement of Dalla Riva, the author and Musolino \cite[Thm.~7.16]{DaLaMu21}. 
\begin{theorem}\label{thm:volmirandao-}
 Let $n\in{\mathbb{N}}\setminus\{0\}$, $\alpha\in]0,1]$. Let $\Omega$ be a bounded open subset of ${\mathbb{R}}^{n}$ of class $C^{1,\alpha}$. Let $r\in]0,+\infty[$ be such that $\overline{\Omega}\subseteq {\mathbb{B}}_n(0,r)$. 
 Let 
 \begin{eqnarray}\label{thm:volmirandao-0a}
{\mathcal{P}}_\Omega^+[k,\varphi] (x)&\equiv& \int_\Omega k(x-y)\varphi(y)\,dy\qquad\forall x\in \overline{\Omega}\,,
\\ \label{thm:volmirandao-0b}
{\mathcal{P}}_\Omega^-[k,\varphi](x)&\equiv& \int_\Omega k(x-y)\varphi(y)\,dy\qquad\forall x\in 
{\mathbb{R}}^n\setminus\Omega\,,
\end{eqnarray}
for all $(k,\varphi)\in {\mathcal{K}}^{1,1}_{-(n-1);o}\times C^{0,\alpha}(\overline{\Omega})$. Then the following statements hold. 
 \begin{enumerate}
\item[(i)] If $\alpha\in]0,1[$, then the bilinear map ${\mathcal{P}}_\Omega^+[\cdot,\cdot]$  from ${\mathcal{K}}^{1,1}_{-(n-1);o}\times C^{0,\alpha}(\overline{\Omega})$ to $C^{1,\alpha}(\overline{\Omega})$ that is delivered by formula (\ref{thm:volmirandao-0a})
is continuous.
\item[(ii)] If $\alpha=1$,  then the bilinear map ${\mathcal{P}}_\Omega^+[\cdot,\cdot]$  from ${\mathcal{K}}^{1,1}_{-(n-1);o}\times C^{0,\alpha}(\overline{\Omega})$ to $C^{1,\omega_1(\cdot)}(\overline{\Omega})$ that is delivered by formula (\ref{thm:volmirandao-0a})
is continuous.

\item[(iii)]  If $\alpha\in]0,1[$, then the bilinear map  ${\mathcal{P}}_\Omega^-[\cdot,\cdot]_{|\overline{{\mathbb{B}}_n(0,r)}\setminus\Omega}$  from ${\mathcal{K}}^{1,1}_{-(n-1);o}\times C^{0,\alpha}(\overline{\Omega})$ to $C^{1,\alpha}(\overline{{\mathbb{B}}_n(0,r)}\setminus\Omega)$ that is delivered by formula (\ref{thm:volmirandao-0b})
is continuous.

\item[(iv)] If $\alpha=1$, then the bilinear map ${\mathcal{P}}_\Omega^-[\cdot,\cdot]_{|\overline{{\mathbb{B}}_n(0,r)}\setminus\Omega}$  from ${\mathcal{K}}^{1,1}_{-(n-1);o}\times C^{0,\alpha}(\overline{\Omega})$ to $C^{1,\omega_1(\cdot)}(\overline{{\mathbb{B}}_n(0,r)}\setminus\Omega)$ that is delivered by formula (\ref{thm:volmirandao-0b})
is continuous.
\end{enumerate}
 \end{theorem}
{\bf Proof.} Since 
\begin{eqnarray*}
\lefteqn{
c'_{\Omega,n-1}\equiv\sup_{x\in\overline{{\mathbb{B}}_n(0,r)}}\int_{
\Omega 
}\frac{dy}{|x-y|^{n-1}} 
}
\\ \nonumber
&&\qquad
\leq 
\sup_{x\in \overline{{\mathbb{B}}_n(0,r)}}
\int_{
{\mathbb{B}}_n(0,r)
}\frac{dy}{|x-y|^{n-1}}
\leq\int_{
{\mathbb{B}}_n(0,2r)
}\frac{dy}{|y|^{n-1}}<+\infty\,,
\end{eqnarray*}
we have
\[
\left|
\int_\Omega k(x-y)\varphi(y)\,dy
\right|\leq c'_{\Omega,n-1}\sup_{\Omega}|\varphi|\sup_{\partial{\mathbb{B}}_n(0,1)}|k|
\qquad\forall x\in \overline{{\mathbb{B}}_n(0,r)} 
\]
for all $k\in {\mathcal{K}}^{1,1}_{-(n-1);o}$ 
 and accordingly
\begin{eqnarray}\label{thm:volmirandao-1}
{\mathcal{P}}_\Omega^+[\cdot,\cdot] 
&\in& {\mathcal{L}}^{(2)}\left(
{\mathcal{K}}^{1,1}_{-(n-1);o}\times C^{0,\alpha}(\overline{\Omega})
,
C^{0}(\overline{\Omega})
\right)\,,
\\ \nonumber
{\mathcal{P}}_\Omega^-[\cdot,\cdot]_{|\overline{{\mathbb{B}}_n(0,r)}\setminus\Omega}
&\in& {\mathcal{L}}^{(2)}\left(
{\mathcal{K}}^{1,1}_{-(n-1);o}\times C^{0,\alpha}(\overline{\Omega})
,
C^{0}(\overline{{\mathbb{B}}_n(0,r)}\setminus\Omega)
\right)\,,
\end{eqnarray}
for all $\alpha\in]0,1]$. 
Let ${\mathbb{E}}$ be a linear and continuous extension map from $C^{0,\alpha}(\overline{\Omega})$ to $C^{0,\alpha}(\overline{{\mathbb{B}}_n(0,r)})$ (cf.~\textit{e.g.}, \cite[Thm.~2.72]{DaLaMu21}).
Then we have
\begin{equation}\label{thm:volmirandao-2a}
 \frac{\partial}{\partial x_l}{\mathcal{P}}_\Omega^+[k,\varphi](x)
=G_l[k,{\mathbb{E}}[\varphi]](x)-{\mathbb{E}}[\varphi] (x) K^+[k,(\nu_\Omega)_l](x)
\end{equation}
for all $x\in   \Omega$ and $(k,\varphi)\in {\mathcal{K}}^{1,1}_{-(n-1);o}\times C^{0,\alpha}(\overline{\Omega})$ (cf.~Dalla Riva, the author and Musolino \cite[Prop.~7.14 (iv)]{DaLaMu21}, where $G_l$ is as in Proposition \ref{prop:gjhc} and
\[
K^+[k,(\nu_\Omega)_l](x)\equiv\int_{\partial\Omega} k(x-y)(\nu_\Omega)_l(y)\,d\sigma_y
\qquad\forall x\in  \Omega \,.
\]
Since  $k\in{\mathcal{K}}^{1,1}_{-(n-1);o}$ and $(\nu_\Omega)_l\in C^{0,\alpha}(\partial\Omega)$, 
an extension of a theorem of Miranda \cite{Mi65}
ensures that $K^+[k,(\nu_\Omega)_l]$ can be extended to a $\alpha$-H\"{o}lder continuous function on $ \overline{\Omega}$ and that $K^+[\cdot,(\nu_\Omega)_l]$ is linear and continuous from ${\mathcal{K}}^{1,1}_{-(n-1);o}$ to $C^{0,\alpha}(\overline{\Omega})$ in case $\alpha\in]0,1[$ (see \cite[Thm.~4.17 (i)]{DaLaMu21}) and   that $K^+[\cdot,(\nu_\Omega)_l]$ is linear and continuous from ${\mathcal{K}}^{1,1}_{-(n-1);o}$ to $C^{0,\omega_1(\cdot)}(\overline{\Omega})$ in case $\alpha=1$ (cf. Theorem \ref{thm:mirandao01} (i) of the Appendix).
Since the pointwise product is 
bilinear and continuous in $C^{0,\alpha}(\overline{\Omega})$
in case $\alpha\in]0,1[$ of statement (i) and in $C^{0,\omega_1(\cdot)}(\overline{\Omega})$ in case $\alpha=1$ of statement (ii),  
the map from ${\mathcal{K}}^{1,1}_{-(n-1);o}\times C^{0,\alpha}(\overline{\Omega})$ to $C^{0,\alpha}(\overline{\Omega})$ in case $\alpha\in]0,1[$ of statement (i) and to $C^{0,\omega_1(\cdot)}(\overline{\Omega})$ in case $\alpha=1$ of statement (ii),  that takes $(k,\varphi)$ to 
${\mathbb{E}}[\varphi]_{|\overline{\Omega}}K^+[\varphi,(\nu_\Omega)_l]$ is bilinear and continuous.
Next we note that the differentiability theorem for integrals depending on a parameter and the Divergence Theorem  imply that
\begin{eqnarray}\label{thm:volmirandao-2b}
\lefteqn{
\frac{\partial}{\partial x_l}{\mathcal{P}}_\Omega^-[k,\varphi](x)
={\mathcal{P}}_\Omega^-[\frac{\partial}{\partial x_l}k,\varphi](x)
}
\\ \nonumber
&&
=
\int_\Omega\frac{\partial}{\partial x_l}k (x-y)(\varphi(y)-{\mathbb{E}}[\varphi](x))\,dy
+{\mathbb{E}}[\varphi](x)\int_\Omega\frac{\partial}{\partial x_l}k(x-y) \,dy
\\ \nonumber
&&
=
\int_\Omega\frac{\partial}{\partial x_l}k (x-y)({\mathbb{E}}[\varphi](y)-{\mathbb{E}}[\varphi](x))\,dy
-{\mathbb{E}}[\varphi](x)\int_{\partial\Omega} k (x-y)(\nu_\Omega)_l(y) \,dy 
\\ \nonumber
&&
=G_l[k,{\mathbb{E}}[\varphi]](x)-{\mathbb{E}}[\varphi] (x) K^-[k,(\nu_\Omega)_l](x)
\end{eqnarray}
for all $x\in  {\mathbb{B}}_n(0,r) \setminus\overline{\Omega}$ and $l\in\{1,\dots,n\}$ and $(k,\varphi)\in {\mathcal{K}}^{1,1}_{-(n-1);o}\times C^{0,\alpha}(\overline{\Omega})$, where $G_l$ is as in Proposition \ref{prop:gjhc} and
\[
K^-[k,(\nu_\Omega)_l](x)\equiv\int_{\partial\Omega} k(x-y)(\nu_\Omega)_l(y)\,d\sigma_y
\qquad\forall x\in {\mathbb{R}}^n\setminus\overline{\Omega} \,.
\]
Since  $k\in{\mathcal{K}}^{1,1}_{-(n-1);o}$ and $(\nu_\Omega)_l\in C^{0,\alpha}(\partial\Omega)$, an extension of a known result of Miranda \cite{Mi65} ensures that $K^-[k,(\nu_\Omega)_l]_{|
{\mathbb{B}}_n(0,r) \setminus\overline{\Omega}
}$ can be extended to a $\alpha$-H\"{o}lder continuous function on $ \overline{{\mathbb{B}}_n(0,r) }\setminus \Omega$ and that $K^-[\cdot,(\nu_\Omega)_l]_{|\overline{{\mathbb{B}}_n(0,r) }\setminus \Omega}$ is linear and continuous from ${\mathcal{K}}^{1,1}_{-(n-1);o}$ to $C^{0,\alpha}(\overline{{\mathbb{B}}_n(0,r)}\setminus\Omega)$ in case $\alpha\in]0,1[$ (see \cite[Thm.~4.17 (ii)]{DaLaMu21}) and that $K^-[k,(\nu_\Omega)_l]_{| {\mathbb{B}}_n(0,r)  \setminus \overline{\Omega}}$ can be extended to a $\omega_1(\cdot)$-H\"{o}lder continuous function on $ \overline{{\mathbb{B}}_n(0,r) }\setminus \Omega$ and that $K^-[\cdot,(\nu_\Omega)_l]_{|\overline{{\mathbb{B}}_n(0,r) }\setminus \Omega}$ is linear and continuous from ${\mathcal{K}}^{1,1}_{-(n-1);o}$ to $C^{0,\omega_1(\cdot)}(\overline{{\mathbb{B}}_n(0,r)}\setminus\Omega)$ in case $\alpha=1$ (cf. Theorem \ref{thm:mirandao01} (ii) of the Appendix). 
Since the pointwise product is 
bilinear and continuous in $C^{0,\alpha}(\overline{{\mathbb{B}}_n(0,r)}\setminus\Omega)$
in case $\alpha\in]0,1[$ of statement (iii) and in $C^{0,\omega_1(\cdot)}(\overline{{\mathbb{B}}_n(0,r)}\setminus\Omega)$ in case $\alpha=1$ of statement (iv),  
the map from ${\mathcal{K}}^{1,1}_{-(n-1);o}\times C^{0,\alpha}(\overline{\Omega})$ to $C^{0,\alpha}(\overline{{\mathbb{B}}_n(0,r)}\setminus\Omega)$ in case $\alpha\in]0,1[$ of statement (iii) and to $C^{0,\omega_1(\cdot)}(\overline{{\mathbb{B}}_n(0,r)}\setminus\Omega)$ in case $\alpha=1$ of statement (iv),  that takes $(k,\varphi)$ to 
${\mathbb{E}}[\varphi]_{|\overline{{\mathbb{B}}_n(0,r)}\setminus\Omega}K[\varphi,(\nu_\Omega)_l]$ is bilinear and continuous.

 By Proposition \ref{prop:gjhc}, $G_l$ is bilinear and continuous. Then the continuity of ${\mathbb{E}}$ and equalities 
(\ref{thm:volmirandao-2a}), (\ref{thm:volmirandao-2b}) imply that 
\begin{eqnarray}\label{thm:volmirandao-3}
&&\frac{\partial}{\partial x_l}{\mathcal{P}}_\Omega^+[k,\varphi] 
\in {\mathcal{L}}^{(2)}\left(
{\mathcal{K}}^{1,1}_{-(n-1);o}\times C^{0,\alpha}(\overline{\Omega})
,
C^{0,\alpha}(\overline{\Omega})
\right)\,,
\\ \nonumber
&&\frac{\partial}{\partial x_l}{\mathcal{P}}_\Omega^-[k,\varphi]_{|\overline{{\mathbb{B}}_n(0,r)}\setminus\Omega}
\in {\mathcal{L}}^{(2)}\left(
{\mathcal{K}}^{1,1}_{-(n-1);o}\times C^{0,\alpha}(\overline{\Omega})
,
C^{0,\alpha}(\overline{{\mathbb{B}}_n(0,r)}\setminus\Omega)
\right)\,,
\end{eqnarray} 
in case  $\alpha\in]0,1[$ of statements (i), (iii) and
\begin{eqnarray}\label{thm:volmirandao-4}
&&\frac{\partial}{\partial x_l}{\mathcal{P}}_\Omega^+[k,\varphi] 
\in {\mathcal{L}}^{(2)}\left(
{\mathcal{K}}^{1,1}_{-(n-1);o}\times C^{0,\alpha}(\overline{\Omega})
,
C^{0,\omega_1(\cdot)}(\overline{\Omega})
\right)\,,
\\ \nonumber
&&\frac{\partial}{\partial x_l}{\mathcal{P}}_\Omega^-[k,\varphi]_{|\overline{{\mathbb{B}}_n(0,r)}\setminus\Omega}
\in {\mathcal{L}}^{(2)}\left(
{\mathcal{K}}^{1,1}_{-(n-1);o}\times C^{0,\alpha}(\overline{\Omega})
,
C^{0,\omega_1(\cdot)}(\overline{{\mathbb{B}}_n(0,r)}\setminus\Omega)
\right)\,,
\end{eqnarray} 
in case  $\alpha=1$ of statements (ii), (iv). Then by combining (\ref{thm:volmirandao-1}) and (\ref{thm:volmirandao-3}), we deduce the validity of  statements (i), (iii) and by combining (\ref{thm:volmirandao-1}) and (\ref{thm:volmirandao-4}), we deduce the validity of  statements (ii), (iv).\hfill  $\Box$ 

\vspace{\baselineskip}

\section{The volume potential with density in a Schauder space of positive exponent}
\label{sec:voscpoe}

We now turn to consider the case in which the density  $\mu$ of the volume potential equals a distribution that is associated to a function of $ L^\infty(\Omega)$, \textit{i.e.}, $\mu={\mathcal{J}}[f]$ with $f\in  L^\infty(\Omega)$
(in the sense of Lemma \ref{lem:cainclo} with any choice of  $\alpha\in]0,1]$). Namely, the so-called classical case, and we first introduce the following classical result. For the convenience of the reader, we include a proof.
\begin{theorem}\label{thm:vopoad}
 Let ${\mathbf{a}}$ be as in (\ref{introd0}), (\ref{ellip}), (\ref{symr}).   Let $S_{ {\mathbf{a}} }$ be a fundamental solution of $P[{\mathbf{a}},D]$. Let $\Omega$ be a bounded open subset of ${\mathbb{R}}^{n}$. The restriction ${\mathcal{P}}_\Omega[S_{ {\mathbf{a}} },\cdot]_{|\overline{ {\mathbb{B}}_n(0,r) }}$ of the  volume potential is linear and continuous from $L^\infty (\Omega)$ to $C^{1,\omega_1(\cdot)}(\overline{ {\mathbb{B}}_n(0,r)  })$ for all $r\in]0,+\infty[$ such that $\overline{\Omega}\subseteq {\mathbb{B}}_n(0,r)$ and formula 
 \begin{equation}\label{thm:vopoad1} 
\frac{\partial}{\partial x_j}{\mathcal{P}}_\Omega[S_{ {\mathbf{a}} },f](x)
=\int_\Omega \frac{\partial S_{ {\mathbf{a}} }}{\partial x_j} (x-y)f(y)\,dy
\qquad\forall x\in {\mathbb{B}}_n(0,r)
\end{equation} holds true for all $f\in L^\infty (\Omega)$, $j\in\{1,\dots,n\}$.
\end{theorem}
{\bf Proof.} Let $r$ be as in the statement. Since $S_{ {\mathbf{a}} }\in C^\infty({\mathbb{R}}^n\setminus\{0\})$ and  both $S_{ {\mathbf{a}} }$ and its first order partial derivatives have a weak singularity in $0$ (cf.~\textit{e.g.} \cite[Lem.~4.2 (i), 4.3 (ii)]{DoLa17}), 
 a classical result implies that the restriction ${\mathcal{P}}_{ {\mathbb{B}}_n(0,r) }[S_{ {\mathbf{a}} },\cdot]_{|\overline{ {\mathbb{B}}_n(0,r) }}$ of the  volume potential ${\mathcal{P}}_{ {\mathbb{B}}_n(0,r) }[S_{ {\mathbf{a}} },\cdot]$ is linear and continuous from $L^\infty ({\mathbb{B}}_n(0,r))$ to $C^1(\overline{ {\mathbb{B}}_n(0,r)  })$ and that formula (\ref{thm:vopoad1}) holds true with $\Omega={\mathbb{B}}_n(0,r)$ for all $f\in L^\infty ({\mathbb{B}}_n(0,r))$ and 
  $x\in {\mathbb{B}}_n(0,r)$  (cf.~\textit{e.g.}, \cite[Prop.~7.5]{DaLaMu21}). Since
 \[
 {\mathcal{P}}_\Omega[S_{ {\mathbf{a}} },f](x)=
 \int_\Omega S_{ {\mathbf{a}} }(x-y)f(y)\,dy= 
 \int_{{\mathbb{B}}_n(0,r)} S_{ {\mathbf{a}} }(x-y)f_\Omega(y)\,dy\qquad\forall
  x\in \Omega\,,
 \]
 for all $f\in L^\infty (\Omega)$, 
where
\begin{equation}\label{thm:vopoad2}
f_\Omega(y)\equiv\left\{
\begin{array}{ll}
 f(y) &\text{if} \ y\in\Omega
\\
0 &\text{if} \ y\in{\mathbb{B}}_n(0,r)\setminus\Omega\,,
\end{array}\right.
\end{equation}
and the map from $L^\infty (\Omega)$ to $L^\infty ({\mathbb{B}}_n(0,r))$ that takes $f$ to $f_\Omega$ is an isometry, we deduce that ${\mathcal{P}}_\Omega[S_{ {\mathbf{a}} },\cdot]_{|\overline{ {\mathbb{B}}_n(0,r) }}$   is linear and continuous from $L^\infty (\Omega)$ to $C^{1}(\overline{ {\mathbb{B}}_n(0,r)  })$  and that formula (\ref{thm:vopoad1}) holds true. Thus it suffices to show that
\begin{equation}\label{thm:vopoad3}
{\mathcal{P}}_\Omega[\frac{\partial S_{ {\mathbf{a}} }}{\partial x_j},\cdot]_{|\overline{ {\mathbb{B}}_n(0,r) }}
\in {\mathcal{L}}\left(L^\infty (\Omega),C^{0,\omega_1(\cdot)}(\overline{ {\mathbb{B}}_n(0,r)  })\right)\qquad\forall j\in \{1,\dots, n\}\,.
\end{equation}
To do so, we wish to apply the abstract result of \cite[Prop.~5.2]{La22a} and we set
\[
X\equiv \overline{ {\mathbb{B}}_n(0,r)  }\,,\quad Y\equiv \Omega\,.
\]
Since
\begin{eqnarray*}
\lefteqn{m_n(
({\mathbb{B}}_n(x,\rho_2)\setminus {\mathbb{B}}_n(x,\rho_1))\cap \Omega
)\leq m_n({\mathbb{B}}_n(0,1))(\rho_2^n-\rho_1^n)
}
\\ \nonumber
&&\qquad\qquad\qquad\qquad\qquad\qquad \forall x\in X,  \rho_1, \rho_2\in [0,+\infty[\ \text{with}\  \rho_1<\rho_2\,,
\end{eqnarray*}
$Y$ is strongly upper $n$-Ahlfors regular with respect to  $X$  in the sense of  \cite[(1.5)]{La22a}.  By \cite[Lem.~4.3]{DoLa17}, we have
\[
\frac{\partial S_{ {\mathbf{a}} }}{\partial x_j}(x-y)\in{\mathcal{K}}_{n-1,n,1}(X\times Y)
\qquad\forall j\in\{1,\dots,n\}\,.
\]
  Next we set 
\[
\upsilon_Y\equiv n\,,\quad
s_1\equiv n-1\,,\quad
s_2\equiv n\,,\quad s_3\equiv 1
\]
  and we note that
  \begin{eqnarray*}
&&s_1\in [\upsilon_Y-1,\upsilon_Y[\,,\ \ 
s_1\geq 0  \,,\ \ 
s_2\in[0,+\infty[\,,\ \  s_3\in]0,1]
\,,
\\
&&C^{0,\max\{r^{\upsilon_Y-s_1},\omega_{s_3}(r)\}}_b(X)=C^{0,\omega_{1}(r)}(X)
\,.
  \end{eqnarray*}
  Then \cite[Prop.~5.2]{La22a} implies that the membership of (\ref{thm:vopoad3}) holds true.\hfill  $\Box$ 

\vspace{\baselineskip}

We are now ready to prove the following extension of a result of Miranda \cite[Thm.~3.I, p.~320]{Mi65} in case $m=0$.
\begin{theorem}\label{thm:vopoadmaz} 
Let ${\mathbf{a}}$ be as in (\ref{introd0}), (\ref{ellip}), (\ref{symr}).   Let $S_{ {\mathbf{a}} }$ be a fundamental solution of $P[{\mathbf{a}},D]$. Let $\alpha\in]0,1]$. 
 Let $\Omega$ be a bounded open subset of  ${\mathbb{R}}^n$ of class $C^{1,\alpha}$. Then the following statements hold. 
  \begin{enumerate}
\item[(i)]  If $\alpha\in ]0,1[$, then
   ${\mathcal{P}}_\Omega^+[S_{ {\mathbf{a}} },\cdot]$ is linear and continuous from $C^{0,\alpha}(\overline{\Omega})$ to $C^{2,\alpha}(\overline{\Omega})$.
\item[(ii)] If $\alpha=1$, then ${\mathcal{P}}_\Omega^+[S_{ {\mathbf{a}} },\cdot]$ is linear and continuous from $C^{0,\alpha}(\overline{\Omega})$ to $C^{2,\omega_1(\cdot)}(\overline{\Omega})$.
\item[(iii)] If $\alpha\in ]0,1[$, then 
 ${\mathcal{P}}_\Omega^-[S_{ {\mathbf{a}} },\cdot]_{|\overline{{\mathbb{B}}_n(0,r)}\setminus\Omega}$ is linear and continuous from $C^{0,\alpha}(\overline{\Omega})$ to the space  $C^{2,\alpha}(\overline{{\mathbb{B}}_n(0,r)}\setminus\Omega)$ for all $r\in]0,+\infty[$ such that $\overline{\Omega}\subseteq {\mathbb{B}}_n(0,r)$.
\item[(iv)] If $\alpha=1$, then ${\mathcal{P}}_\Omega^-[S_{ {\mathbf{a}} },\cdot]_{|\overline{{\mathbb{B}}_n(0,r)}\setminus\Omega}$ is linear and continuous from $C^{0,\alpha}(\overline{\Omega})$ to the space  $C^{2,\omega_1(\cdot)}(\overline{{\mathbb{B}}_n(0,r)}\setminus\Omega)$ for all $r\in]0,+\infty[$ such that $\overline{\Omega}\subseteq {\mathbb{B}}_n(0,r)$.
\end{enumerate} 
\end{theorem}
{\bf Proof.} Let $r\in]0,+\infty[$ be such that $\overline{\Omega}\subseteq {\mathbb{B}}_n(0,r)$. Since 
$C^{0,\alpha}(\overline{\Omega})$ is continuously embedded into $L^\infty(\Omega)$, then Theorem \ref{thm:vopoad} implies that 
\begin{eqnarray}\label{thm:vopoadmaz1} 
&&{\mathcal{P}}_\Omega^+[S_{ {\mathbf{a}} },\cdot]
\in
{\mathcal{L}}\left(C^{0,\alpha}(\overline{\Omega}),C^{0}(\overline{\Omega})\right)\,,
\\ \nonumber
&&{\mathcal{P}}_\Omega^-[S_{ {\mathbf{a}} },\cdot]_{|\overline{{\mathbb{B}}_n(0,r)}\setminus\Omega}
\in
{\mathcal{L}}\left(C^{0,\alpha}(\overline{\Omega}),C^{0}(\overline{{\mathbb{B}}_n(0,r)}\setminus\Omega)
\right) 
\end{eqnarray}
and that
 \begin{eqnarray}\label{thm:vopoadmaz2} 
&&\frac{\partial}{\partial x_j}{\mathcal{P}}_\Omega^+[S_{ {\mathbf{a}} },f](x)
={\mathcal{P}}_\Omega^+[\frac{\partial}{\partial x_j}S_{ {\mathbf{a}} },f](x)
\quad\forall x\in \Omega\,,
\\ \nonumber
&&
\frac{\partial}{\partial x_j}{\mathcal{P}}_\Omega^-[S_{ {\mathbf{a}} },f](x)
={\mathcal{P}}_\Omega^-[\frac{\partial}{\partial x_j}S_{ {\mathbf{a}} },f](x)
\quad\forall x\in {\mathbb{B}}_n(0,r)\setminus\overline{\Omega}
\,,
\end{eqnarray}
for all $f\in C^{0,\alpha}(\overline{\Omega})$, $j\in\{1,\dots,n\}$ both in case $\alpha\in]0,1[$ of statements (i), (iii) and in case $\alpha=1$  of statements (ii), (iv). Thus by the definition of norm in the
spaces 
$C^{2,\alpha}(\overline{\Omega})$ and $C^{2,\alpha}(\overline{{\mathbb{B}}_n(0,r)}\setminus\Omega)$, it suffices to show that
\begin{eqnarray}\label{thm:vopoadmaz3a} 
&&{\mathcal{P}}_\Omega^+[\frac{\partial}{\partial x_j}S_{ {\mathbf{a}} },\cdot]
\in
{\mathcal{L}}\left(C^{0,\alpha}(\overline{\Omega}),C^{1,\alpha}(\overline{\Omega})\right)\,,
\\ \nonumber
&&{\mathcal{P}}_\Omega^-[\frac{\partial}{\partial x_j}S_{ {\mathbf{a}} },\cdot]_{|\overline{{\mathbb{B}}_n(0,r)}\setminus\Omega}
\in
{\mathcal{L}}\left(C^{0,\alpha}(\overline{\Omega}),C^{1,\alpha}(\overline{{\mathbb{B}}_n(0,r)}\setminus\Omega)
\right)\,,
\end{eqnarray}
for all $j\in\{1,\dots,n\}$ in case   $\alpha\in]0,1[$ of statements (i), (iii) and 
\begin{eqnarray}\label{thm:vopoadmaz3b} 
&&{\mathcal{P}}_\Omega^+[\frac{\partial}{\partial x_j}S_{ {\mathbf{a}} },\cdot]
\in
{\mathcal{L}}\left(C^{0,\alpha}(\overline{\Omega}),C^{1,\omega_1(\cdot)}(\overline{\Omega})\right)\,,
\\ \nonumber
&&{\mathcal{P}}_\Omega^-[\frac{\partial}{\partial x_j}S_{ {\mathbf{a}} },\cdot]_{|\overline{{\mathbb{B}}_n(0,r)}\setminus\Omega}
\in
{\mathcal{L}}\left(C^{0,\alpha}(\overline{\Omega}),C^{1,\omega_1(\cdot)}(\overline{{\mathbb{B}}_n(0,r)}\setminus\Omega)
\right)\,,
\end{eqnarray}
for all $j\in\{1,\dots,n\}$ in case   $\alpha=1$ of statements (ii), (iv). Let $j\in\{1,\dots,n\}$. By formula (\ref{prop:grafun1}), we have
\begin{equation}\label{thm:vopoadmaz3aa} 
\frac{\partial}{\partial x_j}S_{ {\mathbf{a}} }(x)=k_{j,1}(x)+k_{j,2}(x)\qquad\forall x\in
{\mathbb{R}}^n\setminus\{0\}\,,
\end{equation}
where
\[
k_{j,1}(x)\equiv \frac{1}{ s_{n}\sqrt{\det a^{(2)} } }
|T^{-1}x|^{-n}\left(x^{t}(a^{(2)})^{-1} \right)_j\qquad\forall x\in
{\mathbb{R}}^n\setminus\{0\} 
\]
and
\[
k_{j,2}(x)\equiv |x|^{2-n}A_{2,j}(\frac{x}{|x|},|x|)+\frac{\partial B_{1}}{\partial x_j}(x)\ln |x|+\frac{\partial C}{\partial x_j} (x)\qquad\forall x\in
{\mathbb{R}}^n\setminus\{0\} \,.
\]
Since $k_{j,1}$ belongs to $C^\infty({\mathbb{R}}^n\setminus\{0\})$ and is positively homogeneous of degree $-(n-1)$ and odd, then
\begin{eqnarray}\label{thm:vopoadmaz4a} 
&&{\mathcal{P}}_\Omega^+[k_{j,1},\cdot]
\in
{\mathcal{L}}\left(C^{0,\alpha}(\overline{\Omega}),C^{1,\alpha}(\overline{\Omega})\right)\,,
\\ \nonumber
&&
{\mathcal{P}}_\Omega^-[k_{j,1},\cdot]_{|\overline{{\mathbb{B}}_n(0,r)}\setminus\Omega}
\in
{\mathcal{L}}\left(C^{0,\alpha}(\overline{\Omega}),C^{1,\alpha}(\overline{{\mathbb{B}}_n(0,r)}\setminus\Omega)
\right)\,,
\end{eqnarray}
in case $\alpha\in]0,1[$ of statements (i), (iii) and
\begin{eqnarray}\label{thm:vopoadmaz4b} 
&&{\mathcal{P}}_\Omega^+[k_{j,1},\cdot]
\in
{\mathcal{L}}\left(C^{0,\alpha}(\overline{\Omega}),C^{1,\omega_1(\cdot)}(\overline{\Omega})\right)\,,
\\ \nonumber
&&
{\mathcal{P}}_\Omega^-[k_{j,1},\cdot]_{|\overline{{\mathbb{B}}_n(0,r)}\setminus\Omega}
\in
{\mathcal{L}}\left(C^{0,\alpha}(\overline{\Omega}),C^{1,\omega_1(\cdot)}(\overline{{\mathbb{B}}_n(0,r)}\setminus\Omega)
\right)\,,
\end{eqnarray}
in case $\alpha=1$  of statements (ii), (iv) 
(cf.  Theorem \ref{thm:volmirandao-}). We now turn to consider the kernel $k_{j,2}$. Since $k_{j,2}$ is of class $C^\infty({\mathbb{R}}^n\setminus\{0\})$, 
\begin{eqnarray*}
\lefteqn{
\sup_{x\in {\mathbb{B}}_n(0,2r)\setminus\{0\}}|x|^{(n-2)+\frac{1}{2}}|k_{j,2}(x)|
\leq \sup_{x\in {\mathbb{B}}_n(0,2r)\setminus\{0\}}
\biggl\{|x|^{\frac{1}{2}}
|A_{2}(\frac{x}{|x|},|x|)| 
}
\\ \nonumber
&&\qquad\qquad\qquad
 +|DB_{1}(x)||x|^{(n-2)+\frac{1}{2}}\ln |x|+|x|^{(n-2)+\frac{1}{2}}|DC(x)|
\biggr\}<+\infty \,,
\end{eqnarray*} 
\[
\sup_{x\in {\mathbb{B}}_n(0,2r)\setminus\{0\}}|x|^{(n-2)+\frac{1}{2}+1}\left|\frac{\partial k_{j,2}}{\partial x_l}(x)\right|<+\infty \qquad\forall l\in\{1,\dots,n\}\,,
\]
and $(n-2)+\frac{1}{2}<n-1$ (cf. Proposition \ref{prop:sedesa}), then a classical result implies that the restriction ${\mathcal{P}}_{ {\mathbb{B}}_n(0,r) }[k_{j,2},\cdot]_{|\overline{ {\mathbb{B}}_n(0,r) }}$ of the  volume potential ${\mathcal{P}}_{ {\mathbb{B}}_n(0,r) }[k_{j,2},\cdot]$ is linear and continuous from $L^\infty ({\mathbb{B}}_n(0,r))$ to $C^1(\overline{ {\mathbb{B}}_n(0,r)  })$ and that the formula
\begin{equation}\label{thm:vopoadmaz5} 
\frac{\partial}{\partial x_l}{\mathcal{P}}_{{\mathbb{B}}_n(0,r)}[k_{j,2},f](x)
=\int_{{\mathbb{B}}_n(0,r)} \frac{\partial k_{j,2}}{\partial x_l} (x-y)f(y)\,dy
\qquad\forall x\in {\mathbb{B}}_n(0,r)
\end{equation} 
holds true for all $f\in L^\infty ({\mathbb{B}}_n(0,r))$, $l\in\{1,\dots,n\}$ (cf.~\textit{e.g.}, \cite[Prop.~7.5]{DaLaMu21}). Since
 \[
 {\mathcal{P}}_\Omega[k_{j,2},f](x)=
 \int_\Omega k_{j,2}(x-y)f(y)\,dy= 
 \int_{{\mathbb{B}}_n(0,r)} k_{j,2}(x-y)f_\Omega(y)\,dy\quad\forall
  x\in {\mathbb{B}}_n(0,r)\,,
 \]
for all $f\in L^\infty (\Omega)$ (see (\ref{thm:vopoad2}) for the definition of $f_\Omega$), and the map from $L^\infty (\Omega)$ to $L^\infty ({\mathbb{B}}_n(0,r))$ that takes $f$ to $f_\Omega$ is an isometry, and $C^{0,\alpha}(\overline{\Omega})$ is continuously embedded into $L^\infty(\Omega)$, we conclude that
\begin{equation}\label{thm:vopoadmaz5a}
{\mathcal{P}}_\Omega[k_{j,2},\cdot]
\in
{\mathcal{L}}\left(C^{0,\alpha}(\overline{\Omega}),C^{1}(\overline{{\mathbb{B}}_n(0,r)})\right)
\end{equation}
and accordingly that
\begin{eqnarray}\label{thm:vopoadmaz6} 
&&{\mathcal{P}}_\Omega^+[k_{j,2},\cdot]
\in
{\mathcal{L}}\left(C^{0,\alpha}(\overline{\Omega}),C^{1}(\overline{\Omega})\right)\,,
\\ \nonumber
&&
{\mathcal{P}}_\Omega^-[k_{j,2},\cdot]_{|\overline{{\mathbb{B}}_n(0,r)}\setminus\Omega}
\in
{\mathcal{L}}\left(C^{0,\alpha}(\overline{\Omega}),C^{1}(\overline{{\mathbb{B}}_n(0,r)}\setminus\Omega)
\right) 
\end{eqnarray}
and
\begin{eqnarray}\label{thm:vopoadmaz7} 
&&\frac{\partial}{\partial x_l}{\mathcal{P}}_\Omega^+[k_{j,2},f](x)
={\mathcal{P}}_\Omega^+[\frac{\partial}{\partial x_l}k_{j,2},f](x)
\quad\forall x\in \Omega\,,
\\ \nonumber
&&
\frac{\partial}{\partial x_l}{\mathcal{P}}_\Omega^-[k_{j,2},f](x)
={\mathcal{P}}_\Omega^-[\frac{\partial}{\partial x_l}k_{j,2},f](x)
\quad\forall x\in {\mathbb{B}}_n(0,r)\setminus\overline{\Omega}
\,,
\end{eqnarray}
for all $\alpha\in]0,1]$. 
We now consider separately cases $\alpha\in]0,1[$ and case $\alpha=1$. We first consider case $\alpha\in]0,1[$ and we wish to prove that
\begin{eqnarray}\label{thm:vopoadmaz8} 
&&{\mathcal{P}}_\Omega^+[\frac{\partial}{\partial x_l}k_{j,2},\cdot]
\in
{\mathcal{L}}\left(C^{0,\alpha}(\overline{\Omega}),C^{0,\alpha}(\overline{\Omega})\right)\,,
\\ \nonumber
&&
{\mathcal{P}}_\Omega^-[\frac{\partial}{\partial x_l}k_{j,2},\cdot]_{|\overline{{\mathbb{B}}_n(0,r)}\setminus\Omega}
\in
{\mathcal{L}}\left(C^{0,\alpha}(\overline{\Omega}),C^{0,\alpha}(\overline{{\mathbb{B}}_n(0,r)}\setminus\Omega)
\right) \,.
\end{eqnarray}
for all $l\in\{1,\dots,n\}$. To do so, we wish to apply the abstract result of \cite[Prop.~5.2]{La22a}. To do so, we set
\[
X^+\equiv\overline{\Omega}\,,\qquad 
X^-\equiv \overline{{\mathbb{B}}_n(0,r)}\setminus \Omega\,,\qquad 
Y\equiv \Omega
\,.
\]
Since
\[
m_n(\Omega\cap  {\mathbb{B}}_n(0,\rho))\leq m_n({\mathbb{B}}_n(0,1))\rho^n
\qquad\forall x\in X^+\cup X^-,  \rho\in ]0,+\infty[\,,
\]
then  $Y$ is upper $n$-Ahlfors regular with respect to both $X^+$ and $X^-$(cf. \textit{e.g.}, \cite[(1.4)]{La22a}).  By Proposition \ref{prop:sedesa} and by the elementary inclusion of \cite[Lem.~3.1]{La22b}, we have
\[
\frac{\partial}{\partial x_l}k_{j,2}\in{\mathcal{K}}_{n-1,n,1}(X^\pm\times Y)
\subseteq {\mathcal{K}}_{n-1,n-(1-\alpha),1-(1-\alpha)}(X^\pm\times Y)
\,,
\]
for all $l\in\{1,\dots,n\}$. Next we set 
\[
\upsilon_Y\equiv n\,,\quad
s_1\equiv n-1\,,\quad
s_2\equiv n-(1-\alpha)\,,\quad s_3\equiv \alpha
\]
  and we note that
  \[
s_1\in [\upsilon_Y-1,\upsilon_Y[\,,\ \ 
s_1\geq 0  \,,\ \ 
s_2\in[0,+\infty[\,,\ \  s_3\in]0,1]
\,,\ \ 
\min\{\upsilon_Y-s_1,s_3\}=\alpha
\,.
  \]
Then \cite[Prop.~5.2]{La22a} implies that
\begin{eqnarray}\label{thm:vopoadmaz9} 
&&{\mathcal{P}}_\Omega^+[\frac{\partial}{\partial x_l}k_{j,2},\cdot]
\in
{\mathcal{L}}\left(L^\infty(\Omega),C^{0,\alpha}(\overline{\Omega})\right)\,,
\\ \nonumber
&&
{\mathcal{P}}_\Omega^-[\frac{\partial}{\partial x_l}k_{j,2},\cdot]_{|\overline{{\mathbb{B}}_n(0,r)}\setminus\Omega}
\in
{\mathcal{L}}\left(L^\infty(\Omega),C^{0,\alpha}(\overline{{\mathbb{B}}_n(0,r)}\setminus\Omega)
\right) \,,
\end{eqnarray}
for all $l\in\{1,\dots,n\}$.  Since $C^{0,\alpha}(\overline{\Omega})$ is continuously embedded into $L^\infty(\Omega)$, the continuity properties of (\ref{thm:vopoadmaz9}) imply the validity of (\ref{thm:vopoadmaz8}). By equality (\ref{thm:vopoadmaz3aa}) and by the memberships of
(\ref{thm:vopoadmaz4a}),  (\ref{thm:vopoadmaz6}), (\ref{thm:vopoadmaz8}), 
we conclude that the memberships of (\ref{thm:vopoadmaz3a}) hold  true and thus the proof of statements (i), (ii) is complete.
We now consider case $\alpha=1$. We wish to prove that
\begin{eqnarray}\label{thm:vopoadmaz8b} 
&&{\mathcal{P}}_\Omega^+[\frac{\partial}{\partial x_l}k_{j,2},\cdot]
\in
{\mathcal{L}}\left(C^{0,\alpha}(\overline{\Omega}),C^{0,\omega_1(\cdot)}(\overline{\Omega})\right)\,,
\\ \nonumber
&&
{\mathcal{P}}_\Omega^-[\frac{\partial}{\partial x_l}k_{j,2},\cdot]_{|\overline{{\mathbb{B}}_n(0,r)}\setminus\Omega}
\in
{\mathcal{L}}\left(C^{0,\alpha}(\overline{\Omega}),C^{0,\omega_1(\cdot)}(\overline{{\mathbb{B}}_n(0,r)}\setminus\Omega)
\right) \,.
\end{eqnarray}
for all $l\in\{1,\dots,n\}$. To do so, we wish to apply the abstract result of \cite[Prop.~5.2]{La22a} and we note that
\begin{eqnarray*}
\lefteqn{m_n(
({\mathbb{B}}_n(x,\rho_2)\setminus {\mathbb{B}}_n(x,\rho_1))\cap \Omega
)\leq m_n({\mathbb{B}}_n(0,1))(\rho_2^n-\rho_1^n)
}
\\ \nonumber
&&\qquad\qquad\qquad\qquad\qquad\qquad \forall x\in X^+\cup X^-,  \rho_1, \rho_2\in [0,+\infty[\ \text{with}\  \rho_1<\rho_2\,.
\end{eqnarray*}
Hence,  $Y$ is strongly upper $n$-Ahlfors regular with respect to both $X^+$ and $X^-$  in the sense of  \cite[(1.5)]{La22a}.  By Proposition \ref{prop:sedesa}, we have
\[
\frac{\partial}{\partial x_l}k_{j,2}\in{\mathcal{K}}_{n-1,n,1}(X^\pm\times Y)
\qquad\forall l\in\{1,\dots,n\}\,.
\]
  Next we set 
\[
\upsilon_Y\equiv n\,,\quad
s_1\equiv n-1\,,\quad
s_2\equiv n\,,\quad s_3\equiv 1
\]
  and we note that
  \begin{eqnarray*}
&&s_1\in [\upsilon_Y-1,\upsilon_Y[\,,\ \ 
s_1\geq 0  \,,\ \ 
s_2\in[0,+\infty[\,,\ \  s_3\in]0,1]
\,,
\\
&&C^{0,\max\{r^{\upsilon_Y-s_1},\omega_{s_3}(r)\}}_b(X^\pm)=C^{0,\omega_{1}(r)}(X^\pm)
\,.
  \end{eqnarray*}
  Then \cite[Prop.~5.2]{La22a} implies that
\begin{eqnarray}\label{thm:vopoadmaz9b} 
&&{\mathcal{P}}_\Omega^+[\frac{\partial}{\partial x_l}k_{j,2},\cdot]
\in
{\mathcal{L}}\left(L^\infty(\Omega),C^{0,\omega_{1}(\cdot)}(\overline{\Omega})\right)\,,
\\ \nonumber
&&
{\mathcal{P}}_\Omega^-[\frac{\partial}{\partial x_l}k_{j,2},\cdot]_{|\overline{{\mathbb{B}}_n(0,r)}\setminus\Omega}
\in
{\mathcal{L}}\left(L^\infty(\Omega),C^{0,\omega_{1}(\cdot)}(\overline{{\mathbb{B}}_n(0,r)}\setminus\Omega)
\right) \,,
\end{eqnarray}
for all $l\in\{1,\dots,n\}$. Since $C^{0,\alpha}(\overline{\Omega})$ is continuously embedded into $L^\infty(\Omega)$, the continuity properties of (\ref{thm:vopoadmaz9b}) imply the validity of (\ref{thm:vopoadmaz8b}). By equality (\ref{thm:vopoadmaz3aa}) and by the memberships of
 (\ref{thm:vopoadmaz4b}), (\ref{thm:vopoadmaz6}), (\ref{thm:vopoadmaz8b}), 
we conclude that the memberships of (\ref{thm:vopoadmaz3b}) hold  true and thus the proof of statements (ii), (iv) is complete.\hfill  $\Box$ 

\vspace{\baselineskip}

Next we introduce the following (known) definition that we need below.
 \begin{definition}\label{defn:vphi} 
  Let ${\mathbf{a}}$ be as in (\ref{introd0}), (\ref{ellip}), (\ref{symr}).   Let $S_{ {\mathbf{a}} }$ be a fundamental solution of $P[{\mathbf{a}},D]$.
Let $\Omega$ be a bounded open  Lipschitz subset of ${\mathbb{R}}^{n}$. If $\phi\in C^{0}(\partial\Omega)$, then we denote by $v_{\Omega}[S_{ {\mathbf{a}} },\phi]$  the single (or simple)  layer potential with moment (or density) $\phi$, i.e., the function from ${\mathbb{R}}^{n}$ to ${\mathbb{R}}$ defined by 
\begin{equation}\label{defn:vphi1} 
v_{\Omega}[S_{ {\mathbf{a}} },\phi](x)\equiv\int_{\partial\Omega}S_{ {\mathbf{a}} }(x-y)\phi (y)\,d\sigma_{y}\qquad\forall x\in {\mathbb{R}}^{n}\,.
\end{equation}
\end{definition}
Under the assumptions of Definition \ref{defn:vphi}, it is known that $v_{\Omega}[S_{ {\mathbf{a}} },\phi]$ is continuous in ${\mathbb{R}}^n$. Indeed,  $\partial\Omega$ is upper $(n-1)$-Ahlfors regular with respect to ${\mathbb{R}}^n$ and $S_{ {\mathbf{a}} }$ has a weak singularity (cf.~\cite[Prop.~6.5]{La24}, \cite[Prop.~4.3]{La22b}, \cite[Lem.~4.2 (i)]{DoLa17}). Then we set
\begin{equation}\label{defn:vphi2}
v_{\Omega}^+[S_{ {\mathbf{a}} },\phi]=v_{\Omega}[S_{ {\mathbf{a}} },\phi]_{|\Omega}\,,\qquad
v_{\Omega}^-[S_{ {\mathbf{a}} },\phi]=v_{\Omega}[S_{ {\mathbf{a}} },\phi]_{|\Omega^-}\,.
\end{equation}

Next, we are ready to prove by induction the following extension of a result of Miranda \cite[Thm.~3.I, p.~320]{Mi65}.
\begin{theorem}\label{thm:vopoadma} 
Let ${\mathbf{a}}$ be as in (\ref{introd0}), (\ref{ellip}), (\ref{symr}).   Let $S_{ {\mathbf{a}} }$ be a fundamental solution of $P[{\mathbf{a}},D]$. 
 Let $m\in {\mathbb{N}}$, $\alpha\in]0,1]$. Let $\Omega$ be a bounded open subset of  ${\mathbb{R}}^n$ of class $C^{m+1,\alpha}$. Then the following statements hold. 
 \begin{enumerate}
 \item[(i)] If $\alpha\in ]0,1[$, then ${\mathcal{P}}_\Omega^+[S_{ {\mathbf{a}} },\cdot]$ is linear and continuous from $C^{m,\alpha}(\overline{\Omega})$ to the space $C^{m+2,\alpha}(\overline{\Omega})$.
 
 \item[(ii)] If $\alpha=1$, then ${\mathcal{P}}_\Omega^+[S_{ {\mathbf{a}} },\cdot]$ is linear and continuous from $C^{m,\alpha}(\overline{\Omega})$ to the space $C^{m+2,\omega_1(\cdot)}(\overline{\Omega})$.
 
\item[(iii)] If $\alpha\in ]0,1[$, then  ${\mathcal{P}}_\Omega^-[S_{ {\mathbf{a}} },\cdot]_{|\overline{{\mathbb{B}}_n(0,r)}\setminus\Omega}$ is linear and continuous from the space $C^{m,\alpha}(\overline{\Omega})$ to the space  $C^{m+2,\alpha}(\overline{{\mathbb{B}}_n(0,r)}\setminus\Omega)$ for all $r\in]0,+\infty[$ such that $\overline{\Omega}\subseteq {\mathbb{B}}_n(0,r)$.

\item[(iv)] If $\alpha=1$, then  ${\mathcal{P}}_\Omega^-[S_{ {\mathbf{a}} },\cdot]_{|\overline{{\mathbb{B}}_n(0,r)}\setminus\Omega}$ is linear and continuous from the space $C^{m,\alpha}(\overline{\Omega})$ to the space  $C^{m+2,\omega_1(\cdot)}(\overline{{\mathbb{B}}_n(0,r)}\setminus\Omega)$ for all $r\in]0,+\infty[$ such that $\overline{\Omega}\subseteq {\mathbb{B}}_n(0,r)$.
\end{enumerate}
\end{theorem}
{\bf Proof.} We proceed by induction on $m$. Case $m=0$ follows by Theorem \ref{thm:vopoadmaz}. We now assume that the statements hold  for $m\geq 0$ and we prove them for $m+1$. By the continuity of the embedding of $C^{m+1,\alpha}(\overline{\Omega})$ into $C^{0,\alpha}(\overline{\Omega})$ and by case $m=0$, we have 
\begin{eqnarray}\label{thm:vopoadma1}
&&{\mathcal{P}}_\Omega^+[S_{ {\mathbf{a}} },\cdot]
\in
{\mathcal{L}}\left(C^{m+1,\alpha}(\overline{\Omega}),C^{0}(\overline{\Omega})\right)\,,
\\ \nonumber
&&{\mathcal{P}}_\Omega^-[S_{ {\mathbf{a}} },\cdot]_{|\overline{{\mathbb{B}}_n(0,r)}\setminus\Omega}\in
{\mathcal{L}}\left(C^{m+1,\alpha}(\overline{\Omega}),C^{0}(\overline{{\mathbb{B}}_n(0,r)}\setminus\Omega)
\right)\,,
\end{eqnarray}
both in case $\alpha\in]0,1[$ of statements (i), (iii) and in case $\alpha=1$  of statements (ii), (iv). 
Since
\[
\sup_{ 0<|\xi|\leq{\mathrm{diam}}\,({\mathbb{B}}_n(0,2r)) }|\xi|^{n-1-1/2}| S_{ {\mathbf{a}} } (\xi)|<+\infty
\]
(cf.~\textit{e.g.}, \cite[Lem.~4.2 (i)]{DoLa17}), 
 formula (\ref{thm:vopoad1}) for the first order derivatives of the volume potential and the integration by parts formula of Theorem \ref{thm:partswsd}  of the Appendix imply that
\begin{eqnarray}\nonumber
\lefteqn{
\frac{\partial}{\partial x_j}{\mathcal{P}}_\Omega[S_{ {\mathbf{a}} },\varphi](x)
=\int_\Omega \frac{\partial S_{ {\mathbf{a}} }}{\partial x_j} (x-y)\varphi (y)\,dy
=-\int_\Omega \frac{\partial S_{ {\mathbf{a}} }}{\partial y_j} (x-y)\varphi (y)\,dy
}
 \\ \label{thm:vopoadma2} 
&&\qquad
=\int_\Omega S_{ {\mathbf{a}} }(x-y)\frac{\partial \varphi}{\partial y_j}\,dy
-\int_{\partial\Omega}S_{ {\mathbf{a}} }(x-y)\varphi(y)(\nu_\Omega)_j(y)\,d\sigma_y
\\ \nonumber
&&\qquad
={\mathcal{P}}_\Omega[S_{ {\mathbf{a}} },\frac{\partial }{\partial y_j}\varphi](x)
-v_\Omega[S_{ {\mathbf{a}} },  (\nu_\Omega)_j\varphi_{|\partial\Omega}](x)
\qquad\forall x\in \Omega\,,
\end{eqnarray}
for all $\varphi\in C^{m+1,\alpha}(\overline{\Omega})$. Moreover, if $x\in{\mathbb{B}}_n(0,r)\setminus\overline{\Omega}$, then $S_{ {\mathbf{a}} }(x-\cdot)\in C^1(\overline{\Omega})$ and thus the Leibnitz rule and the Divergence Theorem imply the validity of  the same equality of (\ref{thm:vopoadma2}).
By the inductive assumption, we have
\begin{eqnarray}\label{thm:vopoadma3} 
&&{\mathcal{P}}_\Omega^+[S_{ {\mathbf{a}} },\frac{\partial }{\partial y_j}(\cdot)]
\in
{\mathcal{L}}\left(C^{m+1,\alpha}(\overline{\Omega}),C^{m+2,\alpha}(\overline{\Omega})\right)\,,
\\ \nonumber
&&{\mathcal{P}}_\Omega^-[S_{ {\mathbf{a}} },\frac{\partial }{\partial y_j}(\cdot)]_{|\overline{{\mathbb{B}}_n(0,r)}\setminus\Omega}\in
{\mathcal{L}}\left(C^{m+1,\alpha}(\overline{\Omega}),C^{m+2,\alpha}(\overline{{\mathbb{B}}_n(0,r)}\setminus\Omega)
\right)\,.
\end{eqnarray} in case $\alpha\in]0,1[$ and
\begin{eqnarray}\label{thm:vopoadma3b} 
&&{\mathcal{P}}_\Omega^+[S_{ {\mathbf{a}} },\frac{\partial }{\partial y_j}(\cdot)]
\in
{\mathcal{L}}\left(C^{m+1,\alpha}(\overline{\Omega}),C^{m+2,\omega_1(\cdot)}(\overline{\Omega})\right)\,,
\\ \nonumber
&&{\mathcal{P}}_\Omega^-[S_{ {\mathbf{a}} },\frac{\partial }{\partial y_j}(\cdot)]_{|\overline{{\mathbb{B}}_n(0,r)}\setminus\Omega}\in
{\mathcal{L}}\left(C^{m+1,\alpha}(\overline{\Omega}),C^{m+2,\omega_1(\cdot)}(\overline{{\mathbb{B}}_n(0,r)}\setminus\Omega)
\right)\,.
\end{eqnarray} in case $\alpha=1$.

Since $\Omega$ is of class $C^{(m+1)+1,\alpha}$, the components of $\nu_\Omega$ are of class $C^{m+1,\alpha}$ 
and the restriction map $r_{|\partial\Omega}[\cdot]$ is linear and continuous from $C^{m+1,\alpha}(\overline{\Omega})$
to $C^{m+1,\alpha}(\partial \Omega)$  and thus the continuity of the pointwise product in 
$C^{m+1,\alpha}(\partial\Omega)$ and a known result for the single layer potential in Schauder spaces imply that
\begin{eqnarray}\label{thm:vopoadma4} 
&&
v_\Omega[S_{ {\mathbf{a}} },  (\nu_\Omega)_jr_{|\partial\Omega}[\cdot]]_{|\overline{\Omega}}
\in
{\mathcal{L}}\left(C^{m+1,\alpha}(\overline{\Omega}),C^{m+2,\alpha}(\overline{\Omega})\right)\,,
\\ \nonumber
&&v_\Omega[S_{ {\mathbf{a}} },  (\nu_\Omega)_jr_{|\partial\Omega}[\cdot]]_{|
\overline{{\mathbb{B}}_n(0,r)}\setminus\Omega}\in
{\mathcal{L}}\left(C^{m+1,\alpha}(\overline{\Omega}),C^{m+2,\alpha}(\overline{{\mathbb{B}}_n(0,r)}\setminus\Omega)
\right)\,,
\end{eqnarray} in case $\alpha\in]0,1[$ of statements (i), (iii)
 (see \cite[Thm.~7.1]{DoLa17}) and
 \begin{eqnarray}\label{thm:vopoadma4b} 
&&
v_\Omega[S_{ {\mathbf{a}} },  (\nu_\Omega)_jr_{|\partial\Omega}[\cdot]]_{|\overline{\Omega}}
\in
{\mathcal{L}}\left(C^{m+1,\alpha}(\overline{\Omega}),C^{m+2,\omega_1(\cdot)}(\overline{\Omega})\right)\,,
\\ \nonumber
&&v_\Omega[S_{ {\mathbf{a}} },  (\nu_\Omega)_jr_{|\partial\Omega}[\cdot]]_{|
\overline{{\mathbb{B}}_n(0,r)}\setminus\Omega}\in
{\mathcal{L}}\left(C^{m+1,\alpha}(\overline{\Omega}),C^{m+2,\omega_1(\cdot)}(\overline{{\mathbb{B}}_n(0,r)}\setminus\Omega)
\right)\,,
\end{eqnarray} in case $\alpha=1$ of statements (ii), (iv) (see Theorem \ref{thm:slaya=1} of the Appendix).

 Then the memberships  of (\ref{thm:vopoadma1}), 
  equality (\ref{thm:vopoadma2}), the memberships of (\ref{thm:vopoadma3}) and (\ref{thm:vopoadma4}) in case $\alpha\in]0,1[$ and the memberships of 
  (\ref{thm:vopoadma3b}) and (\ref{thm:vopoadma4b}) in case $\alpha=1$
  imply the validity of statements (i)--(iv) for (m+1). Thus the induction principle implies the validity of statements (i)--(iv) for all $m\in{\mathbb{N}}$.\hfill  $\Box$ 

\vspace{\baselineskip}

We also note that the following  embedding lemma follows by the classical Theorems  \ref{thm:vopoad},  \ref{thm:vopoadmaz} (i) for the fundamental  solution $S_n$. 
\begin{lemma}\label{lem:amc0b-1a}
  Let $\Omega$ be a bounded open  subset of ${\mathbb{R}}^{n}$. 
  \begin{enumerate}
\item[(i)] If $\alpha\in]0,1[$ and $\Omega$ is of class $C^{0,1}$,  then $L^\infty ( \Omega)$ is continuously embedded into $C^{-1,\alpha}(\overline{\Omega})$.
\item[(ii)] If $\alpha\in]0,1[$ and $\Omega$ is of class $C^{1,\alpha}$,  then $C^{0,\alpha}( \overline{\Omega})$ is continuously embedded into $C^{-1,1}(\overline{\Omega})$.
\end{enumerate}
\end{lemma}
 {\bf Proof.} We first observe that
\[
u=\Delta \left({\mathcal{P}}^+_\Omega[S_n,u]\right)\qquad\forall
u\in L^\infty ( \Omega)\,,
\]
in the sense of distributions. By Theorem \ref{thm:vopoad}, 
 ${\mathcal{P}}^+_\Omega[S_n,\cdot]$ is linear and continuous from $L^\infty( \Omega)$ to $C^{1,\omega_1(\cdot)}(\overline{\Omega})$. Since   we know that  $C^{1,\omega_1(\cdot)}(\overline{\Omega})$ is   continuously embedded into $C^{1,\alpha}(\overline{\Omega})$ and that  $\Delta$ is linear and continuous from $C^{1,\alpha}(\overline{\Omega})$ to
$C^{-1,\alpha}(\overline{\Omega})$, we conclude that  the statement (i) holds true.

By the classical Theorem \ref{thm:vopoadmaz} (i), ${\mathcal{P}}^+_\Omega[S_n,\cdot]$ is linear and continuous from $C^{0,\alpha}(\overline{\Omega})$ to $C^{2,\alpha}(\overline{\Omega})$. Since   we know that  $C^{2,\alpha}(\overline{\Omega})$  is   continuously embedded into $C^{1,1}(\overline{\Omega})$ and that  $\Delta$ is linear and continuous from $C^{1,1}(\overline{\Omega})$ to
$C^{-1,1}(\overline{\Omega})$, we conclude that  the statement (ii) holds true.
\hfill  $\Box$ 

\vspace{\baselineskip}

\section{The volume potential with density in a Schauder space of negative exponent}
\label{sec:voscnege}

Next, we turn to compute the distributional volume potential for the specific form of $\mu$'s in $\left(C^{1,\alpha}(\overline{\Omega})\right)'$ that are extensions of elements of $C^{-1,\alpha}(\overline{\Omega}) $ in the sense of Proposition \ref{prop:nschext}. 

\begin{proposition}\label{prop:dvpsnec-1a}
 Let ${\mathbf{a}}$ be as in (\ref{introd0}), (\ref{ellip}), (\ref{symr}).   Let $S_{ {\mathbf{a}} }$ be a fundamental solution of $P[{\mathbf{a}},D]$. Let $\alpha\in]0,1]$. Let $\Omega$ be a bounded open Lipschitz subset of ${\mathbb{R}}^{n}$. If  $f=  f_{0}+\sum_{j=1}^{n}\frac{\partial}{\partial x_{j}}f_{j}\in C^{-1,\alpha}(\overline{\Omega}) $, then ${\mathcal{P}}_{\Omega}[S_{ {\mathbf{a}} },E^\sharp[f]]$ is the distribution that is associated to the function
 \begin{eqnarray}\label{prop:dvpsnec-1a1}
 \lefteqn{
\int_\Omega S_{ {\mathbf{a}} }(x-y)f_0(y)\,dy
}
\\ \nonumber
&&\qquad
+\sum_{j=1}^{n} \int_{\partial\Omega}S_{ {\mathbf{a}} }(x-y)(\nu_{\Omega})_{j}(y)f_{j}(y)\,d\sigma_y
+\sum_{j=1}^{n}  \frac{\partial}{\partial x_j}\int_\Omega S_{ {\mathbf{a}} }(x-y)f_j(y)\,dy
\end{eqnarray}
 for almost all $x\in{\mathbb{R}}^n$. 
\end{proposition}
{\bf Proof.} If $v\in {\mathcal{D}}({\mathbb{R}}^n)$, then 
\[
\frac{\partial}{\partial y_j}\int_{  {\mathbb{R}}^n }S_{ {\mathbf{a}} }(x-y)v(x)dx
=
\int_{  {\mathbb{R}}^n }\frac{\partial}{\partial y_j}(S_{ {\mathbf{a}} }(x-y))v(x)dx
=-\int_{  {\mathbb{R}}^n }\frac{\partial}{\partial x_j}S_{ {\mathbf{a}} }(x-y)v(x)dx
\]
for all $x\in {\mathbb{R}}^n$. Indeed both $S_{ {\mathbf{a}} }$ and its first order partial derivatives have a weak singularity (cf.
\cite[Lem.~4.2 (i), Lem.~4.3 (ii)]{DoLa17}, \cite[Prop.~7.2, 7.5]{DaLaMu21}). 
Hence,  Proposition \ref{prop:nschext} and the Fubini Theorem imply that
\begin{eqnarray*}
\lefteqn{
\langle {\mathcal{P}}_{\Omega}[E^\sharp[f]],v\rangle =
\langle (r_{|\overline{\Omega}}^tE^\sharp[f])\ast S_{ {\mathbf{a}} },v\rangle 
}
\\ \nonumber
&&\qquad
=\langle E^\sharp[f](y),r_{|\overline{\Omega}}\langle S_{ {\mathbf{a}} }(\eta),v(y+\eta)\rangle\rangle 
\\ \nonumber
&&\qquad
=\int_\Omega f_0(y)\int_{ {\mathbb{R}}^n}S_{ {\mathbf{a}} }(\eta)v(y+\eta)\,d\eta dy
\\ \nonumber
&&\qquad\quad
+
\sum_{j=1}^n\int_{\partial\Omega}f_j(y)(\nu_\Omega)_j(y)
\int_{ {\mathbb{R}}^n}S_{ {\mathbf{a}} }(\eta)v(y+\eta)\,d\eta d\sigma_y
\\ \nonumber
&&\qquad\quad
-\sum_{j=1}^n\int_\Omega f_j(y)\frac{\partial}{\partial y_j}
\int_{  {\mathbb{R}}^n }S_{ {\mathbf{a}} }(\eta)v(y+\eta)\,d\eta dy
\\ \nonumber
&&\qquad
=\int_\Omega f_0(y)\int_{ {\mathbb{R}}^n}S_{ {\mathbf{a}} }(x-y)v(x)\,dx dy
\\ \nonumber
&&\qquad\quad
+
\sum_{j=1}^n\int_{\partial\Omega}f_j(y)(\nu_\Omega)_j(y)
\int_{ {\mathbb{R}}^n}S_{ {\mathbf{a}} }(x-y)v(x)\,dx d\sigma_y
\\ \nonumber
&&\qquad\quad
-\sum_{j=1}^n\int_\Omega f_j(y)\frac{\partial}{\partial y_j}
\int_{  {\mathbb{R}}^n }S_{ {\mathbf{a}} }(x-y)v(x)dx\, dy
\\ \nonumber
&&\qquad
=\int_{ {\mathbb{R}}^n} \int_\Omega S_{ {\mathbf{a}} }(x-y)f_0(y)\,dy\, v(x)dx  
\\ \nonumber
&&\qquad\quad
+
\sum_{j=1}^n\int_{ {\mathbb{R}}^n}\int_{\partial\Omega}f_j(y)(\nu_\Omega)_j(y)S_{ {\mathbf{a}} }(x-y) d\sigma_y\, v(x)dx 
\\ \nonumber
&&\qquad\quad
-\sum_{j=1}^n  \int_\Omega f_j(y)\int_{  {\mathbb{R}}^n }\frac{\partial}{\partial y_j}
 \left(S_{ {\mathbf{a}} }(x-y)\right)v(x)dx\,dy 
 \\ \nonumber
&&\qquad
=\int_{ {\mathbb{R}}^n} \int_\Omega S_{ {\mathbf{a}} }(x-y)f_0(y)dy \, v(x)dx  
\\ \nonumber
&&\qquad\quad
+
\sum_{j=1}^n\int_{ {\mathbb{R}}^n}\int_{\partial\Omega}f_j(y)(\nu_\Omega)_j(y)S_{ {\mathbf{a}} }(x-y) d\sigma_y\, v(x)dx 
\\ \nonumber
&&\qquad\quad
+\sum_{j=1}^n\int_{  {\mathbb{R}}^n } \int_\Omega\frac{\partial}{\partial x_j}
 S_{ {\mathbf{a}} }(x-y)f_j(y)dy\, v(x)dx  
\end{eqnarray*}
and accordingly, ${\mathcal{P}}_{\Omega}[E^\sharp[f]]$ is the distribution that is associated to the function
 in (\ref{prop:dvpsnec-1a1}).\hfill  $\Box$ 

\vspace{\baselineskip}

  Then we have the following generalization to volume potentials of  nonhomogeneous second order elliptic operators of a known result for the Laplace operator (cf.~\cite[Thm.~3.6 (ii)]{La08a}, Dalla Riva, the author and Musolino~\cite[Thm.~7.19]{DaLaMu21})
    \begin{proposition}\label{prop:dvpsnecr-1a}
  Let $\alpha\in]0,1[$.   Let ${\mathbf{a}}$ be as in (\ref{introd0}), (\ref{ellip}), (\ref{symr}).   Let $S_{ {\mathbf{a}} }$ be a fundamental solution of $P[{\mathbf{a}},D]$. Let $\Omega$ be a bounded open subset of ${\mathbb{R}}^{n}$ of class $C^{1,\alpha}$. Let $r\in ]0,+\infty[$ be such that $\overline{\Omega}\subseteq {\mathbb{B}}_n(0,r)$. Then the following statements hold. 
 \begin{enumerate}
\item[(i)] If  $f=  f_{0}+\sum_{j=1}^{n}\frac{\partial}{\partial x_{j}}f_{j}\in C^{-1,\alpha}(\overline{\Omega}) $, then
 \begin{equation}\label{prop:dvpsnecr-1a2}
{\mathcal{P}}_\Omega^+[S_{ {\mathbf{a}} },E^\sharp[f]]\in C^{1,\alpha}(\overline{\Omega}), \ 
{\mathcal{P}}_\Omega^-[S_{ {\mathbf{a}} },E^\sharp[f]]\in C^{1,\alpha}_{{\mathrm{loc}} }(\overline{\Omega^-})  \end{equation}
and 
\begin{equation}\label{defn:Ppm1}
{\mathcal{P}}_\Omega^+[S_{ {\mathbf{a}} },E^\sharp[f]](x)={\mathcal{P}}_\Omega^-[S_{ {\mathbf{a}} },E^\sharp[f]](x)\qquad\forall x\in\partial\Omega\,.
\end{equation}
Moreover,
\begin{eqnarray}\label{prop:dvpsnecr-1a2a}
&&P[{\mathbf{a}},D] {\mathcal{P}}_\Omega^+[S_{ {\mathbf{a}} },E^\sharp[f]]= f\qquad\textit{in}\ {\mathcal{D}}'(\Omega)\,,
\\ \nonumber
&&P[{\mathbf{a}},D] {\mathcal{P}}_\Omega^-[S_{ {\mathbf{a}} },E^\sharp[f]]= 0\qquad\textit{in}\ {\mathcal{D}}'({\mathbb{R}}^n\setminus\overline{\Omega})
\,.
\end{eqnarray}
\item[(ii)] Then the operator ${\mathcal{P}}_\Omega^+[S_{ {\mathbf{a}} },E^\sharp[\cdot]]$ is linear and continuous from $C^{-1,\alpha}(\overline{\Omega}) $ to $C^{1,\alpha}(\overline{\Omega})$.
\item[(iii)]   Then the operator ${\mathcal{P}}_\Omega^-[S_{ {\mathbf{a}} },E^\sharp[\cdot]]_{|\overline{{\mathbb{B}}_n(0,r)}\setminus\Omega}$ is linear and continuous from the space $C^{-1,\alpha}(\overline{\Omega}) $ to   $C^{1,\alpha}(\overline{{\mathbb{B}}_n(0,r)}\setminus\Omega)$.
\end{enumerate}  
\end{proposition}
{\bf Proof.} The equalities in (\ref{prop:dvpsnecr-1a2a}) follow by equalities  (\ref{prop:nschext2}), (\ref{prop:nschext1}) and  (\ref{prop:dvpsa2}). Then equality (\ref{defn:Ppm1}) follows by formula (\ref{prop:dvpsnec-1a1}) for ${\mathcal{P}}_{\Omega}[S_{ {\mathbf{a}} },E^\sharp[f]]$, by the continuity in ${\mathbb{R}}^n$ of the single layer potential with density in $C^{0,\alpha}(\partial\Omega)$ and by the continuous differentiability in ${\mathbb{R}}^n$ of   volume potentials with density in $C^{0,\alpha}(\overline{\Omega})$ (cf. Theorem \ref{thm:vopoad}).

Next, we prove the memberships of (\ref{prop:dvpsnecr-1a2})    and statements (ii), (iii)  
 by exploiting Lemma \ref{lem:coc-1ao} and a variant of the proof of \cite[Thm.~7.19]{DaLaMu21}. 

To do so, we turn to prove that if $(f_0,f_1,\dots,f_n)\in  (C^{0,\alpha}(\overline{\Omega}))^{n+1}$, then the restriction to $\overline{\Omega}$ and to $\overline{{\mathbb{B}}_n(0,r)}\setminus\Omega$ of the function in  (\ref{prop:dvpsnec-1a1}) that is associated to ${\mathcal{P}}_{\Omega}[S_{ {\mathbf{a}} },E^\sharp[f]]$ define    elements of $C^{1,\alpha}(\overline{\Omega})$ and of 
$C^{1,\alpha}(\overline{{\mathbb{B}}_n(0,r)}\setminus\Omega)$, respectively  and that the maps 
$B_+$ and
$B_-$  from  
$(C^{0,\alpha}(\overline{\Omega}))^{n+1}$ to $C^{1,\alpha}(\overline{\Omega})$
and to 
$C^{1,\alpha}(\overline{{\mathbb{B}}_n(0,r)}\setminus\Omega)$ that take  $(f_0,f_1,\dots,f_n)$ to the restriction to $\overline{\Omega}$ and to $\overline{{\mathbb{B}}_n(0,r)}\setminus\Omega$ of the function 
$B [f_0,f_1,\dots,f_n]$ in  (\ref{prop:dvpsnec-1a1}) are linear and continuous, respectively. Here we note that
\begin{eqnarray*}
B_+[f_0,f_1,\dots,f_n]&=&{\mathcal{P}}_\Omega^+[E^\sharp[\Xi[f_0,f_1,\dots,f_n]]]\,,
\\
B_-[f_0,f_1,\dots,f_n]&=&{\mathcal{P}}_\Omega^-[E^\sharp[\Xi[f_0,f_1,\dots,f_n]]]_{|\overline{{\mathbb{B}}_n(0,r)}\setminus\Omega} 
\,,
\end{eqnarray*}
for all $(f_0,f_1,\dots,f_n)\in (C^{0,\alpha}(\overline{\Omega}))^{n+1}$ (cf.~(\ref{eq:xi0a}) for the definition of $\Xi$). 
For the   continuity of the first and third addendum of  (\ref{prop:dvpsnec-1a1}) from $(C^{0,\alpha}(\overline{\Omega}))^{n+1}$ to $C^{1,\alpha}(\overline{\Omega})$ and to 
$C^{1,\alpha}(\overline{{\mathbb{B}}_n(0,r)}\setminus\Omega)$,   we refer to  Theorem  \ref{thm:vopoadmaz} (i), (iii) in case $m=0$.
 
Since $v_\Omega[\cdot]_{|\overline{\Omega}}$ and $v_\Omega[\cdot]_{|\overline{{\mathbb{B}}_n(0,r)}\setminus\Omega}$
are known to be continuous from $C^{0,\alpha}(\partial\Omega)$ to 
$C^{1,\alpha}(\overline{\Omega})$ and to $C^{1,\alpha}(\overline{{\mathbb{B}}_n(0,r)}\setminus\Omega)$, respectively (cf.~\textit{e.g.}, (\ref{defn:vphi2}), \cite[Thm.~7.1 (i)]{DoLa17}), the membership of $\nu_{\Omega}$ in $\left(C^{0,\alpha}(\partial\Omega)\right)^{n} $ and the continuity of the pointwise product in $C^{0,\alpha}(\partial\Omega)$ imply that also
the second addendum of    (\ref{prop:dvpsnec-1a1})  is linear and continuous from $(C^{0,\alpha}(\overline{\Omega}))^{n+1}$ to $C^{1,\alpha}(\overline{\Omega})$ and to $C^{1,\alpha}(\overline{{\mathbb{B}}_n(0,r)}\setminus\Omega)$, respectively. 
In particular, if $f\in C^{-1,\alpha}(\overline{\Omega})$ and 
$ f_{0}+\sum_{j=1}^{n}\frac{\partial}{\partial x_{j}}f_{j}$, then 
\begin{eqnarray}
\nonumber
\lefteqn{{\mathcal{P}}_\Omega^+[S_{ {\mathbf{a}} },E^\sharp[f]]_{|\overline{\Omega}}
={\mathcal{P}}_\Omega^+[S_{ {\mathbf{a}} },E^\sharp[\Xi[f_0,f_1,\dots,f_n]]]_{|\overline{\Omega} }
\in C^{1,\alpha}(\overline{\Omega})\,,
}
\\
\nonumber
\lefteqn{{\mathcal{P}}_\Omega^-[S_{ {\mathbf{a}} },E^\sharp[f]]_{|\overline{{\mathbb{B}}_n(0,r)}\setminus\Omega}
}
\\
\nonumber
&&\qquad\qquad
={\mathcal{P}}_\Omega^-[S_{ {\mathbf{a}} },E^\sharp[\Xi[f_0,f_1,\dots,f_n]]]_{|\overline{{\mathbb{B}}_n(0,r)}\setminus\Omega}
\in C^{1,\alpha}(\overline{{\mathbb{B}}_n(0,r)}\setminus\Omega)  
\end{eqnarray}
and the memberships of (\ref{prop:dvpsnecr-1a2}) hold true. Then Lemma \ref{lem:coc-1ao} implies that statements (ii), (iii) hold true.\hfill  $\Box$ 

\vspace{\baselineskip}

 In case $\alpha=1$,  we know that the elements of $C^{-1,1}(\overline{\Omega})$ are actually essentially  bounded functions (cf.~Proposition \ref{prop:c-1ainfty}) and that $E^\sharp[f]={\mathcal{J}}[f]$ for all $f\in C^{-1,1}(\overline{\Omega})$ (cf.~Proposition \ref{prop:nschexta=1}). Then the classical Theorem \ref{thm:vopoad} on the volume potential implies the validity of the following statement.
  \begin{proposition}\label{prop:dvpsnecr-11}
 Let ${\mathbf{a}}$ be as in (\ref{introd0}), (\ref{ellip}), (\ref{symr}).   Let $S_{ {\mathbf{a}} }$ be a fundamental solution of $P[{\mathbf{a}},D]$. Let $\Omega$ be a bounded open Lipschitz subset of ${\mathbb{R}}^{n}$. Let $r\in ]0,+\infty[$ be such that $\overline{\Omega}\subseteq {\mathbb{B}}_n(0,r)$. Then the following statements hold. 
 \begin{enumerate}
\item[(i)] If  $f=  f_{0}+\sum_{j=1}^{n}\frac{\partial}{\partial x_{j}}f_{j}\in C^{-1,1}(\overline{\Omega}) $, then
\begin{eqnarray}\label{prop:dvpsnecr-110}
&&{\mathcal{P}}_\Omega^+[S_{ {\mathbf{a}} },E^\sharp[f]]=
{\mathcal{P}}_\Omega^+[S_{ {\mathbf{a}} },{\mathcal{J}}[f]] \,,
\\
\label{prop:dvpsnecr-111}
&&{\mathcal{P}}_\Omega^+[S_{ {\mathbf{a}} },E^\sharp[f]]\in C^{1,\omega_1(\cdot)}(\overline{\Omega}), \ 
{\mathcal{P}}_\Omega^-[S_{ {\mathbf{a}} },E^\sharp[f]]\in C^{1,\omega_1(\cdot)}_{{\mathrm{loc}} }(\overline{\Omega^-})  \end{eqnarray}
and equalities (\ref{defn:Ppm1}) and  (\ref{prop:dvpsnecr-1a2a})
 are satisfied. 
\item[(ii)] The operator ${\mathcal{P}}_\Omega^+[S_{ {\mathbf{a}} },E^\sharp[\cdot]]$ is linear and continuous from $C^{-1,1}(\overline{\Omega}) $ to the space $C^{1,\omega_1(\cdot)}(\overline{\Omega})$.
\item[(iii)] The operator ${\mathcal{P}}_\Omega^-[S_{ {\mathbf{a}} },E^\sharp[\cdot]]_{|\overline{{\mathbb{B}}_n(0,r)}\setminus\Omega}$ is linear and continuous from $C^{-1,1}(\overline{\Omega}) $ to  $C^{1,\omega_1(\cdot)}(\overline{{\mathbb{B}}_n(0,r)}\setminus\Omega)$.
\end{enumerate}  
\end{proposition}
{\bf Proof.} By Proposition \ref{prop:nschexta=1}, equality (\ref{prop:dvpsnecr-110}) holds true. Since  Proposition \ref{prop:c-1ainfty} implies that $C^{-1,1}(\overline{\Omega})$ is continuously embedded into $L^\infty(\Omega)$, then the classical Theorem \ref{thm:vopoad} implies that ${\mathcal{P}}_\Omega[S_{ {\mathbf{a}} },E^\sharp[\cdot]]={\mathcal{P}}_\Omega^+[S_{ {\mathbf{a}} },{\mathcal{J}}[\cdot]]={\mathcal{P}}_\Omega^+[S_{ {\mathbf{a}} },\cdot]$ is linear and continuous from $C^{-1,1}(\overline{\Omega}) $ to $C^{1,\omega_1(\cdot)}(\overline{{\mathbb{B}}_n(0,r)})$. Then the continuity of the restriction operators from $C^{1,\omega_1(\cdot)}(\overline{{\mathbb{B}}_n(0,r)})$ to $C^{1,\omega_1(\cdot)}(\overline{\Omega})$ and  to  $C^{1,\omega_1(\cdot)}(\overline{{\mathbb{B}}_n(0,r)}\setminus\Omega)$ implies the validity of statements (ii), (iii). \hfill  $\Box$ 

\vspace{\baselineskip}


\section*{Declarations}

\begin{itemize}
\item Funding: the author  acknowledges  the support of the Research  Project GNAMPA-INdAM   $\text{CUP}\_$E53C22001930001 `Operatori differenziali e integrali in geometria spettrale' and   of the Project funded by the European Union – Next Generation EU under the National Recovery and Resilience Plan (NRRP), Mission 4 Component 2 Investment 1.1 - Call for tender PRIN 2022 No. 104 of February, 2 2022 of Italian Ministry of University and Research; Project code: 2022SENJZ3 (subject area: PE - Physical Sciences and Engineering) ``Perturbation problems and asymptotics for elliptic differential equations: variational and potential theoretic methods''.	
\item  Conflict of interest/Competing interests: this paper does not have any  conflict of interest or competing interest.  
\item Ethics approval and consent to participate: not applicable.
\item Consent for publication: not applicable.
\item Data availability: not applicable. 
\item Materials availability: not applicable.
\item Code availability: not applicable. 
\item Author contribution: applies to the entire paper.
\end{itemize}

\section{Appendix: a formula of integration by parts}
 
If  $X$ and $Y$ are subsets of ${\mathbb{R}}^n$, $s\in {\mathbb{R}}$, then we denote by ${\mathcal{K}}_{s,X\times Y}$, the set of continuous functions $K$ from $(X\times Y)\setminus {\mathbb{D}}_{ X\times Y }$ to ${\mathbb{C}}$ such that
\[
 \|K\|_{ {\mathcal{K}}_{s,X\times Y} }\equiv \sup_{(x,y)\in  (X\times Y)\setminus {\mathbb{D}}_{ X\times Y }  }|K(x,y)|\,|x-y|^s<+\infty\,.
\]
Then we prove the following formula of integration by parts. Related formulas
 are known. See for example Mitrea, Mitrea and Mitrea \cite[Thm.~1.11.8]{MitMitMit22}. For the convenience of the reader, we include a proof.
 \begin{theorem}\label{thm:partswsd}
Let $n\in{\mathbb{N}}\setminus\{0,1\}$.   Let $\Omega$ be a nonempty bounded open Lipschitz subset of  ${\mathbb{R}}^n$. Let $s\in ]0,n-1]$. Let $K\in {\mathcal{K}}_{s, \Omega\times \overline{\Omega}}$, $\varphi\in C^0(\overline{\Omega})\cap C^1(\Omega)$.  Let $x\in \Omega$. Let $j\in\{1,\dots,n\}$. Let $\frac{\partial K}{\partial y_j}(x,\cdot)$ exist and be continuous in $\overline{\Omega}\setminus\{x\}$. 
Let $\frac{\partial\varphi}{\partial y_j}\in L^1(\Omega)$. If $s=(n-1)$, we further assume that 
\begin{equation}\label{thm:partswsd1}
\Psi_j(K,x)\equiv \lim_{\epsilon\to 0}\int_{\partial{\mathbb{B}}_n(0,1)}K(x,x-\epsilon\xi)\xi_j\,d\sigma_\xi \epsilon^{n-1}
\end{equation}
exists and belongs to ${\mathbb{C}}$. Then the following statements hold.
\begin{enumerate}
\item[(i)] If $s\in ]0,n-1[$, then 
\[
\Psi_j(K,x)\equiv \lim_{\epsilon\to 0}\int_{\partial{\mathbb{B}}_n(0,1)}K(x,x-\epsilon\xi)\xi_j\,d\sigma_\xi \epsilon^{n-1}=0\,.
\]
\item[(ii)] The function $K(x,\cdot)\frac{\partial \varphi }{\partial y_j}(\cdot)$ is integrable in $\Omega$.

\item[(iii)] The principal value
\[
{\mathrm{p.v.}}\int_\Omega\frac{\partial K}{\partial y_j}(x,y)\varphi(y)\,dy
\equiv\lim_{\epsilon\to 0}\int_{\Omega
\setminus {\mathbb{B}}_n(x,\epsilon)
}\frac{\partial K}{\partial y_j}(x,y)\varphi(y)\,dy
\]
 exists in ${\mathbb{C}}$ and 
\begin{eqnarray}\label{thm:partswsd2}
\lefteqn{
{\mathrm{p.v.}}\int_\Omega\frac{\partial K}{\partial y_j}(x,y)\varphi(y)\,dy
}
\\ \nonumber
&&\qquad
=- \int_\Omega K(x,y)\frac{\partial \varphi }{\partial y_j} (y)\,dy+\int_{\partial \Omega}K(x,y)\varphi(y)(\nu_\Omega)_j(y)\,d\sigma_y
\\ \nonumber
&&\qquad\quad
+\varphi(x)\Psi_j(K,x)\,.
\end{eqnarray}
\end{enumerate}
\end{theorem}
{\bf Proof.} (i) It suffices to observe that
\begin{eqnarray*}
\lefteqn{
\left|\int_{\partial{\mathbb{B}}_n(0,1)}K(x,x-\epsilon\xi)\xi_j\,d\sigma_\xi \epsilon^{n-1}
\right|
}
\\ \nonumber
&&\qquad
\leq \|K\|_{{\mathcal{K}}_{s, \Omega\times \overline{\Omega}}}\int_{\partial{\mathbb{B}}_n(0,1)}|\epsilon\xi|^{-s}|\xi_j|\,d\sigma_\xi \epsilon^{n-1}
\leq \|K\|_{{\mathcal{K}}_{s, \Omega\times \overline{\Omega}
}} s_n \epsilon^{(n-1)-s} 
\end{eqnarray*}
for all $\epsilon \in ]0,{\mathrm{dist}}\,(x,\partial \Omega)[$.

(ii) If $\epsilon \in ]0,{\mathrm{dist}}\,(x,\partial \Omega)[$, then $\overline{{\mathbb{B}}_{n}(x,\epsilon)}\subseteq \Omega$ and the set 
\[
\Omega_{\epsilon}\equiv \Omega\setminus
\overline{{\mathbb{B}}_{n}(x,\epsilon) }
\]
is of Lipschitz class. Then we note that the function $K(x,\cdot)\frac{\partial \varphi }{\partial y_j}(\cdot)$ is measurable in $\Omega$ and that
\begin{eqnarray*}
\lefteqn{
\int_{\Omega}\left|K(x,y)\frac{\partial \varphi }{\partial y_j}(y)\right|\,dy
}
\\ \nonumber
&&\qquad
\leq\int_{\Omega_\epsilon} \frac{\|K\|_{{\mathcal{K}}_{s, \Omega\times \overline{\Omega}}}}{\epsilon^{s}}
\left|\frac{\partial \varphi }{\partial y_j}(y)\right|\,dy
+\int_{{\mathbb{B}}_n(x,\epsilon)} \frac{\|K\|_{{\mathcal{K}}_{s, \Omega\times \overline{\Omega}}}}{|x-y|^s}\sup_{{\mathbb{B}}_n(x,\epsilon)}\left|\frac{\partial \varphi }{\partial y_j} \right|\,dy
\\ \nonumber
&&\qquad
\leq \|K\|_{{\mathcal{K}}_{s, \Omega\times \overline{\Omega}}}
\left\{\epsilon^{-s}
\left\|\frac{\partial \varphi }{\partial y_j}\right\|_{L^1(\Omega)}
+\sup_{{\mathbb{B}}_n(x,\epsilon)}\left|\frac{\partial \varphi }{\partial y_j} \right|
\int_{{\mathbb{B}}_n(x,\epsilon)}\frac{dy}{|x-y|^s}
\right\} <+\infty\,.
\end{eqnarray*}
Hence, statement (ii) is true.

(iii) By the Divergence Theorem (cf.~\textit{e.g.}, \cite[Thm.~4.1]{DaLaMu21}), we have
\begin{eqnarray*}
\lefteqn{
\int_{\Omega_\epsilon}\frac{\partial K}{\partial y_j}(x,y)\varphi(y)\,dy
}
\\ \nonumber
&&\qquad
=-\int_{\Omega_\epsilon}K(x,y)\frac{\partial \varphi }{\partial y_j} (y)\,dy
+\int_{\Omega_\epsilon}\frac{\partial  }{\partial y_j}(K(x,y)\varphi(y))\,dy
\\ \nonumber
&&\qquad
=-\int_{\Omega_\epsilon}K(x,y)\frac{\partial \varphi }{\partial y_j} (y)\,dy
+\int_{\partial \Omega_\epsilon}K(x,y)\varphi(y)(\nu_{\Omega_\epsilon})_j(y)\,d\sigma_y
\\ \nonumber
&&\qquad
=-\int_{\Omega_\epsilon}K(x,y)\frac{\partial \varphi }{\partial y_j} (y)\,dy
+\int_{\partial \Omega }K(x,y)\varphi(y)(\nu_{\Omega })_j(y)\,d\sigma_y
\\ \nonumber
&&\qquad\quad
+\int_{\partial {\mathbb{B}}_{n}(x,\epsilon)}K(x,y)\varphi(y)
\frac{x_{j}-y_{j}}{|x-y|}\,d\sigma_{y}
\qquad \forall\epsilon\in]0,{\mathrm{dist}}\,(x,\partial \Omega)[\,.
\end{eqnarray*}
By (ii) and by the Dominated Convergence Theorem, we have
\[
\lim_{\epsilon\to 0}\int_{\Omega_\epsilon}K(x,y)\frac{\partial \varphi }{\partial y_j} (y)\,dy
=\int_{\Omega}K(x,y)\frac{\partial \varphi }{\partial y_j} (y)\,dy\,.
\]
Next we note that
\begin{eqnarray*}
\lefteqn{
\int_{\partial {\mathbb{B}}_{n}(x,\epsilon)}K(x,y)\varphi(y)
\frac{x_{j}-y_{j}}{|x-y|}\,d\sigma_{y}
=\varphi(x)\int_{\partial {\mathbb{B}}_{n}(x,\epsilon)}K(x,y) 
\frac{x_{j}-y_{j}}{|x-y|}\,d\sigma_{y}
}
\\ \nonumber
&&\qquad\qquad\qquad\qquad\qquad
+\int_{\partial {\mathbb{B}}_{n}(x,\epsilon)}K(x,y)(\varphi(y)-\varphi(x))
\frac{x_{j}-y_{j}}{|x-y|}\,d\sigma_{y}\,,
\end{eqnarray*}
that
\begin{eqnarray*}
\lefteqn{
\int_{\partial {\mathbb{B}}_{n}(x,\epsilon)}K(x,y) 
\frac{x_{j}-y_{j}}{|x-y|}\,d\sigma_{y}
}
\\ \nonumber
&&\qquad
=\int_{\partial {\mathbb{B}}_{n}(0,1)}K(x,x-\epsilon\xi) 
\frac{\xi_{j}}{|\xi|}\,d\sigma_{\xi}\epsilon^{n-1}=\int_{\partial{\mathbb{B}}_n(0,1)}K(x,x-\epsilon\xi)\xi_j\,d\sigma_\xi\epsilon^{n-1}\,,
\end{eqnarray*}
and that
\begin{eqnarray*}
\lefteqn{
\left|\int_{\partial {\mathbb{B}}_{n}(x,\epsilon)}K(x,y)(\varphi(y)-\varphi(x))
\frac{x_{j}-y_{j}}{|x-y|}\,d\sigma_{y}\right|
}
\\ \nonumber
&&\qquad
\leq\sup_{y\in \partial {\mathbb{B}}_{n}(x,\epsilon)}|\varphi(y)-\varphi(x)|
\|K\|_{{\mathcal{K}}_{s, \Omega\times \overline{\Omega}}}
\int_{\partial {\mathbb{B}}_{n}(x,\epsilon)}|x-y|^{-s}\,d\sigma_y
\\ \nonumber
&&\qquad
=\sup_{y\in \partial {\mathbb{B}}_{n}(x,\epsilon)}|\varphi(y)-\varphi(x)|
\|K\|_{{\mathcal{K}}_{s, \Omega\times \overline{\Omega}}}s_{n}\epsilon^{(n-1)-s}
\end{eqnarray*}
for all $ \epsilon \in ]0,{\mathrm{dist}}\,(x,\partial \Omega)[$. Then by taking the limit as $\epsilon$ tends to zero
and by the definition of $\Psi_j(K,x)$,  we deduce the validity of the formula of (iii).\hfill  $\Box$ 

\vspace{\baselineskip}

\section{Appendix: a limiting case of two theorems of C.~Miranda}
We now introduce the following extension  to the limiting case in which $\alpha=1$ of a classical result of Miranda \cite{Mi65} (see also \cite[Thm.~4.17]{DaLaMu21}), who has considered the case of domains of class $C^{1,\alpha}$  and of densities $\mu \in C^{0,\alpha}(\partial\Omega)$ for $\alpha\in]0,1[$. For a proof, we refer to Dalla Riva, the author and Musolino \cite{DaLaMu24a}.
\begin{theorem}\label{thm:mirandao01}
Let $\Omega$ be a bounded open subset of ${\mathbb{R}}^n$ of class $C^{1,1}$. Then the following statements hold.
\begin{enumerate}
\item[(i)] For each $(k,\mu)\in {\mathcal{K}}^{1,1}_{-(n-1);o}\times C^{0,1}(\partial\Omega)$, the map
\[
  \int_{\partial\Omega} k(x-y) \mu(y)\,d\sigma_{y}\qquad \forall x\in\Omega 
\]
can be extended to a unique  $\omega_1(\cdot)$-H\"{o}lder continuous function $K[k,\mu]^+$ on $\overline{\Omega}$. Moreover, the map from ${\mathcal{K}}^{1,1}_{-(n-1);o}\times C^{0,1}(\partial\Omega)$ to $C^{0,\omega_1(\cdot)}(\overline{\Omega})$ that  takes $(k,\mu)$ to $K[k,\mu]^+$ is bilinear and continuous.
\item[(ii)] Let $r\in]0,+\infty[$ be such that $\overline{\Omega}\subseteq {\mathbb{B}}_n(0,r)$. Then for each   $(k,\mu)\in {\mathcal{K}}^{1,1}_{-(n-1);o}\times C^{0,1}(\partial\Omega)$  the map
\[
  \int_{\partial\Omega} k(x-y) \mu(y)\,d\sigma_{y}\qquad \forall x\in
 {\mathbb{R}}^n\setminus\overline{\Omega}\,,
\]
can be extended to a unique continuous function $K[k,\mu]^-$ on ${\mathbb{R}}^n\setminus \Omega$ such that 
the restriction $K[k,\mu]^-_{|\overline{{\mathbb{B}}_n(0,r)}\setminus\Omega}$ is $\omega_1(\cdot)$-H\"{o}lder continuous. Moreover,  the map from ${\mathcal{K}}^{1,1}_{-(n-1);o}\times C^{0,1}(\partial\Omega)$ to $C^{0,\omega_1(\cdot)}( \overline{{\mathbb{B}}_n(0,r)}\setminus\Omega)$ that takes $(k,\mu)$ to $K[k,\mu]^-_{|\overline{{\mathbb{B}}_n(0,r)}\setminus\Omega}$ is bilinear and continuous.
\end{enumerate}
\end{theorem}
 Next we consider an  extension to the limiting case $\alpha=1$ of a classical result of Miranda~\cite{Mi65} for the single layer potential (see also  
Wiegner~\cite{Wi93}, Dalla Riva \cite{Da13}, Dalla Riva, Morais and Musolino \cite{DaMoMu13}). For a proof, we refer to \cite{La24e}.
\begin{theorem}\label{thm:slaya=1}
 Let ${\mathbf{a}}$ be as in (\ref{introd0}), (\ref{ellip}), (\ref{symr}).   Let $S_{ {\mathbf{a}} }$ be a fundamental solution of $P[{\mathbf{a}},D]$. 
 Let $m\in {\mathbb{N}}\setminus\{0\}$. Let $\Omega$ be a bounded open subset of  ${\mathbb{R}}^n$ of class $C^{m,1}$.  Then the following statements hold. 
 \begin{enumerate}
 \item[(i)]  If $\mu\in C^{m-1,1}(\partial\Omega)$, then the function 
 $v_\Omega^+[S_{ {\mathbf{a}} },\mu]$  belongs to $C^{m,\omega_1(\cdot)}(\overline{\Omega})$
 and the operator from $C^{m-1,1}(\partial\Omega)$ to $C^{m,\omega_1(\cdot)}(\overline{\Omega})$ that takes $\mu$  to $v_\Omega^+[S_{ {\mathbf{a}} },\mu]$   is linear and continuous.
 \item[(ii)] Let $r\in]0,+\infty[$ be such $\overline{\Omega}\subseteq {\mathbb{B}}_n(0,r)$.
  If $\mu\in C^{m-1,1}(\partial\Omega)$, then the function $v_\Omega^-[S_{ {\mathbf{a}} },\mu]_{|\overline{{\mathbb{B}}_n(0,r)}\setminus\Omega}$ belongs
  to the space  $C^{m,\omega_1(\cdot)}(\overline{{\mathbb{B}}_n(0,r)}\setminus\Omega)$ 
  and the operator from $C^{m-1,1}(\partial\Omega)$ to $C^{m,\omega_1(\cdot)}(\overline{{\mathbb{B}}_n(0,r)}\setminus\Omega)$ that takes $\mu$  to $v_\Omega^-[S_{ {\mathbf{a}} },\mu]_{|\overline{{\mathbb{B}}_n(0,r)}\setminus\Omega}$   is linear and continuous.
 \end{enumerate}
\end{theorem}

\end{document}